\documentclass[10pt]{amsart}
\usepackage[foot]{amsaddr}
\usepackage{graphicx}
\graphicspath{images/}

\usepackage{amsmath}
\usepackage{amsfonts}
\usepackage{tikz-cd} 
\usepackage{amssymb}
\usepackage{amsthm}
\usepackage{tikzsymbols}
\usepackage{mathtools}
\usepackage{xcolor}
\usepackage{adjustbox}
\usepackage{blkarray}
\usepackage{makecell}
\usepackage{float}
\usepackage{caption}
\usepackage[font={small,it}]{caption}
\usepackage[margin=1cm]{caption}
\usepackage{longtable}
\usepackage{tikz-3dplot}
\usepackage{stmaryrd}
\usetikzlibrary{patterns}

\newcommand{\taures}{\tau^{\text{res}}}
\newcommand{\taudel}{\tau^{\text{del}}}
\newcommand{\taucont}{\tau^{\text{cont}}}
\newcommand{\rres}{r^{\text{res}}}
\newcommand{\rdel}{r^{\text{del}}}
\newcommand{\rcont}{r^{\text{cont}}}

\newcommand{\cal}{\mathcal}
\newcommand{\class}[1]{\widetilde{\mathcal{C}}'_{#1}}
\newcommand{\betaInv}{\widetilde{\beta}}
\newcommand{\coloop}{U_{1,1}}
\newcommand{\unif}{U_{a,b}}
\newcommand{\dunif}{U_{b-a,b}}
\newcommand{\nat}{M_\rho^k}

\usepackage{theoremref}
\usepackage{hyperref}

\theoremstyle{plain}
\newtheorem{theorem}{Theorem}[subsection]

\newtheorem{proposition}[theorem]{Proposition}
\newtheorem{lemma}[theorem]{Lemma}
\newtheorem{corollary}[theorem]{Corollary}
\newtheorem{observation}[theorem]{Observation}

\theoremstyle{definition}
\newtheorem{definition}{Definition}[subsection]
\newtheorem{example}[definition]{Example}

\newtheorem{remark}[definition]{Remark}

\usepackage{blindtext}
\usepackage{diagbox}

\usepackage[many]{tcolorbox}
\tcolorboxenvironment{lemma}{
  breakable,
  colback=blue!10!white,
  boxrule=0pt,
  boxsep=1pt,
  left=2pt,right=2pt,top=2pt,bottom=2pt,
  oversize=2pt,
  sharp corners,
  before skip=\topsep,
  after skip=\topsep,
}

\tcolorboxenvironment{observation}{
  breakable,
  colback=blue!10!white,
  boxrule=0pt,
  boxsep=1pt,
  left=2pt,right=2pt,top=2pt,bottom=2pt,
  oversize=2pt,
  sharp corners,
  before skip=\topsep,
  after skip=\topsep,
}

\tcolorboxenvironment{example}{
  breakable,
  colback=green!10!white,
  boxrule=0pt,
  boxsep=1pt,
  left=2pt,right=2pt,top=2pt,bottom=2pt,
  oversize=2pt,
  sharp corners,
  before skip=\topsep,
  after skip=\topsep,
}

\tcolorboxenvironment{strategy}{
  breakable,
  colback=black!10!white,
  boxrule=0pt,
  boxsep=1pt,
  left=2pt,right=2pt,top=2pt,bottom=2pt,
  oversize=2pt,
  sharp corners,
  before skip=\topsep,
  after skip=\topsep,
}

\tcolorboxenvironment{corollary}{
  breakable,
  colback=blue!10!white,
  boxrule=0pt,
  boxsep=1pt,
  left=2pt,right=2pt,top=2pt,bottom=2pt,
  oversize=2pt,
  sharp corners,
  before skip=\topsep,
  after skip=\topsep,
}

\tcolorboxenvironment{theorem}{
  breakable,
  colback=blue!10!white,
  boxrule=0pt,
  boxsep=1pt,
  left=2pt,right=2pt,top=2pt,bottom=2pt,
  oversize=2pt,
  sharp corners,
  before skip=\topsep,
  after skip=\topsep,
}

\tcolorboxenvironment{proposition}{
  breakable,
  colback=blue!10!white,
  boxrule=0pt,
  boxsep=1pt,
  left=2pt,right=2pt,top=2pt,bottom=2pt,
  oversize=2pt,
  sharp corners,
  before skip=\topsep,
  after skip=\topsep,
}

\tcolorboxenvironment{remark}{
  breakable,
  colback=green!10!white,
  boxrule=0pt,
  boxsep=1pt,
  left=2pt,right=2pt,top=2pt,bottom=2pt,
  oversize=2pt,
  sharp corners,
  before skip=\topsep,
  after skip=\topsep,
}

\tcolorboxenvironment{definition}{
  breakable,
  colback=green!10!white,
  boxrule=0pt,
  boxsep=1pt,
  left=2pt,right=2pt,top=2pt,bottom=2pt,
  oversize=2pt,
  sharp corners,
  before skip=\topsep,
  after skip=\topsep,
}

\usepackage{subfiles} 

\title[The essential bound of a polymatroid]{The essential bound of a polymatroid and its applications to excluded minor problems}
\author[F. Young]{Fiona Young}
\address{Department of Mathematics, Cornell University, Ithaca, New York, 14853}
\email{fy79@cornell.edu}
\date{\today}

\begin{document}

\maketitle

ABSTRACT. The singleton and doubleton minors of a polymatroid $\rho$ encode a surprising amount of information about the structural complexity of $\rho$. Given any polymatroid $\rho$, we can subtract from it a maximally-separated polymatroid, resulting in a $k$-polymatroid. We introduce a notion of boundedness for $\rho$ that corresponds to $k$. Our results provide an organized framework for thinking about polymatroid excluded minor problems. In particular, let $\mathcal{C}$ denote the minor-closed class of matroids characterized by excluding the uniform matroid $U_{2,b}$ and its dual $U_{b-2,b}$. We show that the list of excluded minors for the class of $k$-polymatroids whose $k$-natural matroids are in $\mathcal{C}$ is finite. We also investigate the more general case of excluding $U_{a,b}$ and its dual $U_{b-a,b}$.



\tableofcontents

\section{Introduction}
\label{section:Introduction}

Just as there are many different ways to define the matroid, the same is true for its multiset analog, the integer polymatroid (which we abbreviate as \textit{polymatroid}). In this paper we primarily use two cryptomorphic perspectives: (1) polymatroids as rank functions, and (2) polymatroids as polytopes.

The rank function definition is extremely versatile, and it is where we will begin.  A \textit{polymatroid} is a pair $(E,\rho)$ where $E$ is a finite set and $\rho: 2^E \to \mathbb{Z}_{\geq 0}$ is a function satisfying:
\begin{enumerate}
    \item $\rho(\varnothing) = 0$.
    \item \textit{Monotonicity}: If $A \subseteq B \subseteq E$, then $\rho(A) \leq \rho(B)$.
    \item \textit{Submodularity}: For all $A, B \subseteq E$, $\rho(A) + \rho(B) \geq \rho(A \cap B) + \rho(A \cup B)$.
\end{enumerate}
If $\rho(\{e\}) \leq k$ for all $e \in E$, then $\rho$ is a \textit{k-polymatroid}. 

The polytopal view, on the other hand, provides the intuition that is instrumental in leading us to our main results. Let $(E,\rho)$ be a polymatroid. The \textit{base polytope of $\rho$} is the set
    \begin{align*}
        B_\rho := \left\{\textbf{x} \in \mathbb{R}^E: \sum\limits_{e \in E}x_e = \rho(E), \sum\limits_{e \in A}x_e \leq \rho(A) \text{ for all } A \subseteq E\right\},
    \end{align*}
    and the \textit{independent set polytope of $\rho$} is the set
    \begin{align*}
        I_\rho := \left\{\textbf{x} \in \mathbb{R}^E: 0 \leq \sum\limits_{e \in A}x_e \leq \rho(A) \text{ for all } A \subseteq E\right\}.
    \end{align*}

We briefly mention a third perspective of the polymatroid -- its geometric representation. Here, polymatroids generalize matroids by allowing elements of higher rank. While matroids consist only of loops (elements of rank $0$) and points (rank $1$), $k$-polymatroids can also contain lines (rank $2$), planes (rank $3$), etc., up to elements of rank $k$.


To analyze the structure of a polymatroid, we turn to classic matroid theory for inspiration. Many important classes of matroids are \textit{minor-closed}, meaning they are closed under the operations of deletion and contraction. We characterize such a class $\cal{C}$ by its set of \textit{excluded minors}: matroids not in $\cal{C}$ whose proper minors are all in $\cal{C}$. An example is the class of $\mathbb{F}$-representable matroids, for a fixed field $\mathbb{F}$. Rota's conjecture states that if $\mathbb{F}$ is finite, then the set of excluded minors for $\mathbb{F}$-representable matroids is finite. Geelen, Gerards, and Whittle \cite{geelen gerards whittle} announced a proof of this conjecture in 2014.

It is natural to ask the same question for polymatroids. Characterizing the class of $\mathbb{F}$-representable $k$-polymatroids by its excluded minors appears to be a difficult task. For the simplest nontrivial case when $k = 2$ and $\mathbb{F} = \mathbb{F}_2$, i.e., the class of \textit{binary} $2$-polymatroids, Oxley, Semple, and Whittle constructed an infinite family of excluded minors \cite{wheelswhirls}. We consider a variation of this problem by first assigning a unique matroid to each $k$-polymatroid: its \textit{$k$-natural matroid}. Geometrically, the \textit{$k$-natural matroid} of a polymatroid $(E,\rho)$ is obtained by replacing each $e \in E$ with $k$ points lying freely in $e$. A similar notion is the \textit{natural matroid} of $(E,\rho)$, obtained by replacing each $e \in E$ by $\rho(\{e\})$ points lying freely in $E$.

If a class $\cal{C}$ of matroids is minor-closed, then so are the following:
\begin{enumerate}
    \item The class $\cal{C}'$ of polymatroids whose natural matroids are in $\cal{C}$.
    \item The class $\cal{C}'_k$ of $k$-polymatroids in $\cal{C}'$.
    \item The class $\class{k}$ of $k$-polymatroids whose $k$-natural matroids are in $\cal{C}$.
\end{enumerate}

Bonin and Long \cite{bonin long} determined the set of excluded minors for $\cal{C}'_2$, where $\cal{C}$ is the class of binary matroids, characterized by the single identically self-dual excluded minor $U_{2,4}$. They found an infinite sequence of excluded minors, along with eight other excluded minors that do not belong to this sequence. One might ask the same question for $\cal{C}'_k$ and $\class{k}$, with $k \geq 3$. As for the latter, the complete (finite) list of excluded minors can be found in \cite{young}. The polymatroid operation of $l$-compression, a `hybrid' of deletion and contraction, was a crucial tool in these efforts. In particular, it suggests that we look deeper into the connection between the structure of $\rho$ and that of $\nat$. For this task, the polytopal view can be quite helpful -- granting us a new framework with which to tackle excluded minor problems for polymatroids.

In Section \ref{s-{Definitions and background}}, we provide the necessary definitions and background for polymatroids and their $k$-natural matroids as rank functions. In Section \ref{s-{Polymatroids as polytopes}}, we transition to the polytopal perspective, and make precise the connection between the two perspectives. This lays the groundwork for Section \ref{s-{Essential boundedness for polymatroids}}, where we prove the main result of this paper which says, in some sense, that the data of the singleton and doubleton minors of a polymatroid $\rho$ suffice in bounding the total complexity of $\rho$. This result has an immediate application in the aforementioned excluded minor problems, which we discuss in Section \ref{s-{Applications to excluded minor problems}}. In particular, when $\cal{C}$ is the minor-closed class of matroids characterized by $Ex(\cal{C}) = \{U_{2,b}, U_{b-2,b}\}$, our main result implies that $Ex\big(\class{k}\big)$ is finite.

\section*{Acknowledgement}
The author is grateful to Joe Bonin for carefully reading the draft of this paper and for identifying some mistakes and typos.

\section{Definitions and background}
\label{s-{Definitions and background}}

Many matroid notions carry over to polymatroids. We follow Oxley \cite{oxley} for standard matroid terminology and notation. In this section we summarize the polymatroid concepts and definitions that are crucial to the results in this paper. For additional background information on polymatroids, we encourage the reader to refer to \cite{bonin long} and \cite{young}.

A brief comment on the various types of interval notation used in this paper: 
\begin{itemize}
    \item A set of double brackets denotes an interval in $\mathbb{R}$. (The set of all real numbers between $a$ and $b$ inclusive is written as $[[a,b]]$.)
    \item A set of single brackets denotes an interval in $\mathbb{Z}$. (The set of all integers between $a$ and $b$ inclusive is written as $[a,b]$.)
    \item We will not use the common notation $[a]$ which refers to the set of integers $\{1, \hdots, a\}$. This is because we will often need to append $0$ to this set, which we can accomplish using the notation $[0,a]$ introduced above. 
\end{itemize}

\begin{definition}
    Let $(E,\rho)$ be a polymatroid and let $X \subseteq E$. 
    \begin{enumerate}
        \item  The \textit{deletion} of $X$ is the polymatroid $(E-X, \rho_{\backslash X})$ where \\$(\rho_{\backslash X})(Y) = \rho(Y)$ for all $Y \subseteq E-X$. The deletion of $X$ from $\rho$ is equivalent to the \textit{restriction} of $\rho$ to $E-X$, denoted $(E-X, \rho|_{E-X})$.
    \item The \textit{contraction} of $X$ is the polymatroid $(E-X, \rho_{/X})$ where \\$(\rho_{/X})(Y) = \rho(X \cup Y)-\rho(X)$ for all $Y \subseteq E-X$. 
    \end{enumerate}
\end{definition}
Deletion and contraction commute when applied to disjoint subsets. A polymatroid $\rho'$ is a \textit{minor} of $\rho$ if $\rho'$ can be obtained from $\rho$ via a sequence of deletions and contractions on disjoint subsets.

\begin{definition}\label{def-{nullity}}
       The \textit{nullity} of a polymatroid $(E,\rho)$ is
\begin{displaymath}
    \eta(\rho) := |E|-\rho(E).
\end{displaymath}
\end{definition}

\subsection{The \texorpdfstring{$k$}{}-natural matroid of a \texorpdfstring{$k$}{}-polymatroid}
\label{ss-{The \texorpdfstring{$k$}{}-natural matroid of a \texorpdfstring{$k$}{}-polymatroid}}

\begin{definition}\label{def-{k-natural matroid}}
    Let $(E,\rho)$ be a polymatroid. For each $e \in E$, let $X_e$ be a set of $k$ elements. If $e, f \in E$ are distinct, then we require $X_e \cap X_f = \varnothing$. For any subset $A \subseteq E$, we define 
        \begin{displaymath}
            X_A := \bigcup_{e \in A} X_e.
        \end{displaymath}
        The \textit{$k$-natural matroid} $\nat$ of $\rho$ is the matroid $(X_E, r)$ where 
        \begin{displaymath}
            r(X) := \min\{\rho(A) + |X-X_A|: A \subseteq E\}.
        \end{displaymath}
\end{definition}

For the reader who prefers geometric representations: $\nat$ is obtained from $\rho$ by replacing each $e \in E$ with $k$ points lying freely in $e$.

By a result of McDiarmid \cite{mcdiarmid}, we see that $\nat  = (X_E, r)$ is indeed a matroid. In general, we say that two elements $x,y \in X_E$ are \textit{clones} if the transposition that swaps $x$ and $y$ while fixing all other elements of $X_E$ is an automorphism of $\nat$. If every pair of distinct elements $a, b \in X \subseteq E$ are clones, then we say that $X$ is a \textit{set of clones}.


\begin{lemma}[\cite{bonin long}]\label{lem-{k-natural is unique}}
    Let $\rho$, $E$, $X_e$, and $X_A$ be as above. A matroid $M$ on $X_E$ is $\nat$ if and only if each set $X_e$ is a set of clones and $r_M(X_A) = \rho(A)$ for all $A \subseteq E$.
\end{lemma}

\begin{observation}[\cite{bonin long}]\label{obs-{k-natural observations}} \leavevmode
\begin{enumerate}
    \item If $\rho$ and $\rho'$ are polymatroids on $E$ with $\nat  = M_{\rho'}^k$ where for each $e \in E$ the corresponding set $X_e$ is the same in both $k$-natural matroids, then $\rho = \rho'$. 
    \item For any polymatroid $\rho$, we have $M_{\rho\backslash \{e\}}^k = M_{\rho}^k\backslash X_e$ and $M_{\rho/\{e\}}^k = \nat /X_e$ for all $e \in E$.
\end{enumerate}
\end{observation}

\subsection{Compression}
\label{ss-{Compression}}

\begin{definition}
\label{def-{compression}}
    Let $(E,\rho)$ be a polymatroid. For $l \in [0,\rho(\{e\})]$, the \textit{$l$-compression of $\rho$ by $e \in E$} is the polymatroid $\rho_{\downarrow e}^l$ obtained by freely adding $l$ points $e_1,\hdots, e_l$ to $e$, then contracting $e_1, \hdots, e_l$ and deleting $e$. If $\rho' = \rho_{\downarrow e}^l$ for some $e$ and $l$, then we say that $\rho'$ is a \textit{compression} of $\rho$ (and $\rho$ is a \textit{decompression} of $\rho'$). If $l$ is known, we can also say that $\rho'$ is an \textit{$l$-compression} of $\rho$ (and $\rho$ is an \textit{$l$-decompression} of $\rho'$). 
    
    Note: When $l = 0$, $\rho_{\downarrow e}^l = \rho_{\backslash \{e\}}$. Similarly, when $l = \rho(\{e\})$, $\rho_{\downarrow e}^l = \rho_{/\{e\}}$. An \textit{internal compression} of $\rho$ is $\rho_{\downarrow e}^l$ where $l \in [1, \rho(\{e\})-1]$.
\end{definition}

A polymatroid $\rho'$ is a \textit{c-minor} \cite{bonin long} of a polymatroid $\rho$ if $\rho'$ can be obtained from $\rho$ through a sequence of deletions, contractions, and compressions on disjoint subsets.

\begin{lemma}[\cite{bonin long}]\label{lem-{k-natural of l-compression}}
    For a $k$-polymatroid $(E,\rho)$ and $l \in [1,k-1]$, fix $e \in E$ with $\rho(\{e\}) > 0$ and fix $e_1,\hdots,e_l \in X_e$. The $k$-natural matroid of $\rho_{\downarrow e}^l$ is 
    \begin{align*}
        \nat /\{e_1,\hdots,e_l\}\backslash (X_e-\{e_1,\hdots,e_l\}).
    \end{align*}
\end{lemma}

\subsection{Duality}
\label{ss-{Duality}}

There are many ways to extend the notion of duality to polymatroids. As in \cite{bonin long} and \cite{young}, we use $k$-duality, which behaves well with the polymatroid operations of deletion, contraction, and compression. 

\begin{definition}
    For a $k$-polymatroid $(E,\rho)$, its \textit{$k$-dual} is the $k$-polymatroid $(E, \rho^*)$ where
\begin{displaymath}
    \rho^*(X): = k|X| + \rho(E-X)-\rho(E).
\end{displaymath}
\end{definition}

\begin{lemma}[\cite{bonin long}]\label{lem-{k-dual of k-natural}}
    Let $(E,\rho)$ be a $k$-polymatroid and $(E,\rho^*)$ be its $k$-dual. The $k$-natural matroid of $\rho$ is dual to that of $\rho^*$, i.e., $(\nat )^* = M_{\rho^*}^k$.
\end{lemma}

\begin{lemma}[\cite{bonin long}]
    Let $(E,\rho)$ be a $k$-polymatroid and $(E,\rho^*)$ be its $k$-dual. If $e \in E$, then $M_{\rho_{\backslash \{e\}}}^k = \left( M_{\rho^*_{/\{e\}}}^k\right)^*$ and $M_{\rho_{/\{e\}}}^k = \left( M_{\rho^*_{\backslash \{e\}}}^k\right)^*$.
\end{lemma}

\begin{theorem}[\cite{young}]\label{thm-{k-duality}}
Let $\cal{C}$ be a minor-closed, dual-closed matroid class. Then $\widetilde{\cal{C}}'_k$ is closed under $k$-duality, as is its set of excluded minors.
\end{theorem}

\section{Polymatroids as polytopes}
\label{s-{Polymatroids as polytopes}}

We now turn towards the polytopal perspective.

\begin{definition}\label{def-{polymatroid as polytopes}}
    Let $(E,\rho)$ be a polymatroid. The \textit{base polytope of $\rho$} is the set
    \begin{align*}
        B_\rho := \left\{\textbf{x} \in \mathbb{R}^E: \sum\limits_{e \in E}x_e = \rho(E), \sum\limits_{e \in A}x_e \leq \rho(A) \text{ for all } A \subseteq E\right\},
    \end{align*}
    and the \textit{independent set polytope of $\rho$} is the set
    \begin{align*}
        I_\rho := \left\{\textbf{x} \in \mathbb{R}^E: 0 \leq \sum\limits_{e \in A}x_e \leq \rho(A) \text{ for all } A \subseteq E\right\}.
    \end{align*}
\end{definition}

\begin{example}
    Let $(\{e, f, g\},\rho)$ be the $3$-polymatroid such that the rank of each singleton is $3$, the rank of each doubleton is $5$, and the total rank is $6$. Its base polytope $B_\rho$ and independent set polytope $I_\rho$ are shown below.

$$
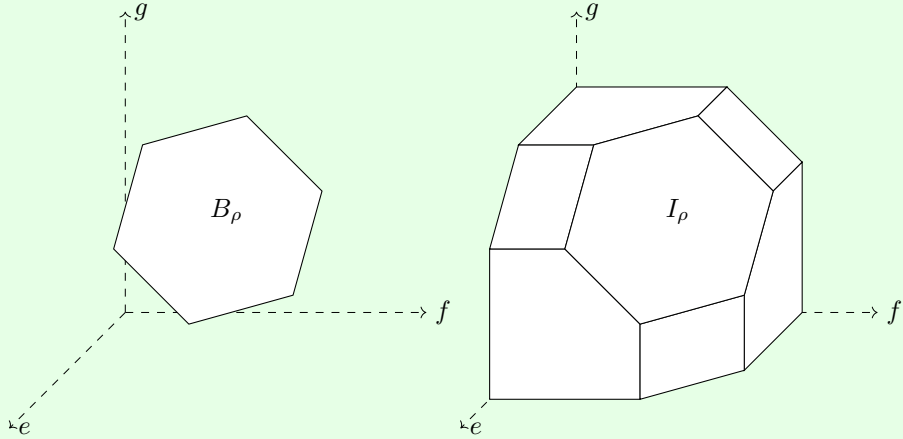
\begin{tikzpicture}
		[ lightdashed/.style={very thin,gray, dashed},
            darkdashed/.style={very thick, black, dashed},
            darksolid/.style = {very thick, black},
			axis/.style={->,black, very thin, dashed}]
			
	\draw[axis] (0,0,0) -- (4,0,0) node[anchor=west]{$f$};
	\draw[axis] (0,0,0) -- (0,4,0) node[anchor=west]{$g$};
	\draw[axis] (0,0,0) -- (0,0,4) node[anchor=west]{$e$};

 \draw[fill = white] (1,2,3)--(1,3,2)--(2,3,1)--(3,2,1)--(3,1,2)--(2,1,3)--cycle;

 \node[] at (2.5,2.5,3) {$B_\rho$}; 
 \end{tikzpicture}\begin{tikzpicture}
		[ lightdashed/.style={very thin,gray, dashed},
            darkdashed/.style={very thick, black, dashed},
            darksolid/.style = {very thick, black},
			axis/.style={->,black, very thin, dashed}]
			
	\draw[axis] (0,0,0) -- (4,0,0) node[anchor=west]{$f$};
	\draw[axis] (0,0,0) -- (0,4,0) node[anchor=west]{$g$};
	\draw[axis] (0,0,0) -- (0,0,4) node[anchor=west]{$e$};

 \draw[fill = white] (0,0,3)--(0,2,3)--(1,2,3)--(2,1,3)--(2,0,3)--cycle;
 \draw[fill = white](0,2,3)--(1,2,3)--(1,3,2)--(0,3,2)--cycle;
 \draw[fill = white] (1,3,2)--(2,3,1)--(2,3,0)--(0,3,0)--(0,3,2)--cycle;
 \draw[fill = white] (2,3,1)--(3,2,1)--(3,2,0)--(2,3,0)--cycle;
 \draw[fill = white] (3,2,1)--(3,1,2)--(3,0,2)--(3,0,0)--(3,2,0)--cycle;
 \draw[fill = white] (3,0,2)--(2,0,3)--(2,1,3)--(3,1,2)--cycle;

 \draw[fill = white] (1,2,3)--(1,3,2)--(2,3,1)--(3,2,1)--(3,1,2)--(2,1,3)--cycle;
 \node[] at (2.5,2.5,3) {$I_\rho$}; 
 \end{tikzpicture}$$

The \textit{permutohedron} $P_n$ is a polytope in $\mathbb{R}^n$ constructed by taking the convex hull of the $n!$ points of the form $(\sigma(1), \hdots, \sigma(n))$ where $\sigma$ is a permutation of $[1,n]$. The polytope $B_\rho$ is equivalent to the permutohedron $P_3$. 
\end{example}

Postnikov \cite{postnikov} introduced \textit{generalized permutohedra} -- deformations of $P_n$ that preserve the directions of its edges, not precluding the possibility that some of its faces may completely degenerate. The set of generalized permutohedra and the set of base polytopes of polymatroids are equivalent up to translation. On the other hand, the independent set polytope is not a generalized permutohedron; it is a \textit{Q-polytope}, which is defined and discussed more in \cite{ardila benedetti doker}.

From the polytopal perspective, we can build new polymatroids out of smaller ones via the operation of Minkowski summation.
\begin{definition}
    The \textit{Minkowski sum} of two polytopes $P$ and $Q$ in $\mathbb{R}^E$ is:
                \begin{align*}
                P+Q = \{\textbf{p}+\textbf{q}: \textbf{p} \in P, \textbf{q} \in Q\}.
                \end{align*}
    Let $R \subseteq \mathbb{R}^E$ be a polytope such that $R+Q = P$. This allows us to define the \textit{Minkowski difference} as $P-Q = R$.
\end{definition}

The set of generalized permutohedra and the set of $Q$-polytopes are each closed under Minkowski summation, as are the set of polymatroid base polytopes and the set of polymatroid independent set polytopes \cite{ardila benedetti doker}.

\subsection{A dissection of \texorpdfstring{$I_\rho$}{}}

By imposing some additional structure on $I_\rho$, we gain an advantage -- the ability to simultaneously track $c$-minors at the level of $\rho$ and at the level of $\nat$.

\begin{definition}\label{def-{multiset rank function}}
Let $(E,\rho)$ be a $k$-polymatroid. Recall that the ground set of $\nat$ is $X_E$ (see Definition \ref{def-{k-natural matroid}} for details). Let us denote the elements of $X_E$ as follows:
\begin{align*}
    X_E &= \{e_i: e \in E, i \in [1,k]\}.
\end{align*} 
    Let $L_\rho^k := [0,k]^E \subseteq \mathbb{R}^E$. Define the following maps:
    \begin{enumerate}
        \item The \textit{partition map}
        \begin{align*}
            \pi_\rho: 2^{X_E} &\to L_\rho^k\\
            X &\mapsto (|X \cap X_e|)_{e \in E}
        \end{align*}
        \item The \textit{multiset rank function} 
        \begin{align*}
            R_\rho: L_\rho^k &\to \mathbb{Z}_{\geq 0}\\
            \textbf{a}&\mapsto \max\limits_{\textbf{b} \in I_\rho \cap \prod\limits_{e \in E}[[0,a_e]]}\left(\sum\limits_{e \in E}b_e\right) .
        \end{align*}
    \end{enumerate}
\end{definition}

\begin{example} \label{ex-{example1}}
Let $\rho$ be the $3$-polymatroid on $\{e,f\}$ such that $\rho(\{e\}) = 3$, $\rho(\{f\}) = 2$, and $\rho(\{e, f\}) = 4$. The ground set of $M_\rho^3$ is
\begin{align*}
    X_E^3 & =  \{e_1,e_2,e_3,f_1,f_2,f_3\}.
\end{align*}
Let $X_1 = \{e_2,f_1,f_2,f_3\}$ and $X_2 = \{e_1,e_3,f_2\}$. We have $\pi_\rho(X_1) = (1,3)$ and $\pi_\rho(X_2) = (2,1)$. The multiset rank function gives $R_\rho(1,3) = R_\rho(2,1) = 3$.
$$
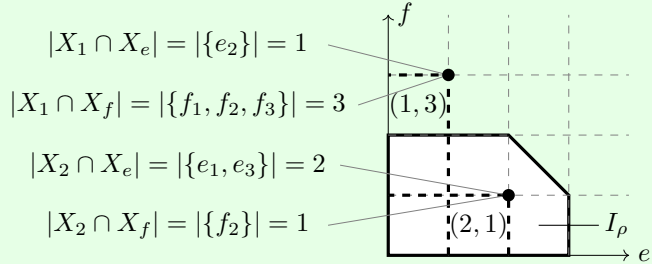
\begin{tikzpicture}
		[scale = 0.8, lightdashed/.style={very thin,gray, dashed},
            darkdashed/.style={very thick, black, dashed},
            darksolid/.style = {very thick, black},
            lightsolid/.style = {very thin, gray},
			axis/.style={->,black, very thin}]
			
	\draw[axis] (0,0)--(4,0) node[anchor=west]{$e$};
	\draw[axis] (0,0) -- (0,4) node[anchor=west]{$f$};

    \draw[darksolid, fill = white] (0,0)--(3,0)--(3,1)--(2,2)--(0,2)--cycle;

    \draw[lightdashed] (1,0)--(1,4);
    \draw[lightdashed] (2,0)--(2,4);
    \draw[lightdashed] (3,0)--(3,4);
    \draw[lightdashed] (0,1)--(4,1);
    \draw[lightdashed] (0,2)--(4,2);
    \draw[lightdashed] (0,3)--(4,3);

    \node[] at (1.5,0.5) {$(2,1)$}; 
    \draw[darkdashed] (2,0)--(2,1)--(0,1);
    \node[] at (0.5,2.5) {$(1,3)$}; 
    \draw[darkdashed] (1,0)--(1,3)--(0,3);

    \node[] at (-3.5,3.5) {$|X_1 \cap X_e| = |\{e_2\}| = 1$}; 
    \draw[lightsolid] (1,3)--(-1,3.5) {};
    \node[] at (-3.5,2.5) {$|X_1 \cap X_f| = |\{f_1,f_2,f_3\}| = 3$}; 
    \draw[lightsolid] (1,3)--(-0.5,2.5) {};
    \node[] at (-3.5,1.5) {$|X_2 \cap X_e| = |\{e_1,e_3\}| = 2$}; 
    \draw[lightsolid] (2,1)--(-0.8,1.5) {};
    \node[] at (-3.5,0.5) {$|X_2 \cap X_f| = |\{f_2\}| = 1$}; 
    \draw[lightsolid] (2,1)--(-1,0.5) {};

    \node[shape = circle, fill = black, scale = 0.5] at (2,1)  {}; 
    \node[shape = circle, fill = black, scale = 0.5] at (1,3)  {}; 

    \draw[] (2.5,0.5)--(3.5,0.5) {};
    \node[] at (3.8,0.5) {$I_\rho$};

 \end{tikzpicture}$$
\end{example}

The next theorem explicitly describes the relationship between our construction above and the $k$-natural matroid of $\rho$ as it was defined in Section \ref{s-{Definitions and background}}.
\begin{theorem}\label{thm-{two ways to define the natural matroid}}
    The $k$-natural matroid $\nat$ of $\rho$ is equal to $(X_E, R_\rho \circ \pi_\rho)$. Equivalently, the following diagram commutes.
    \begin{displaymath}
       \begin{tikzcd}
2^{X_E} \arrow[rd, "\pi_\rho"'] \arrow[rr, "r_{\nat }"] && \mathbb{Z}_{\geq 0} \\
& {L_\rho^k} \arrow[ru, "R_\rho"'] &                 
\end{tikzcd}
    \end{displaymath}
\end{theorem}
\begin{proof}
Let $X \subseteq X_E$. Abbreviate $I_\rho^X : = I_\rho \cap \prod\limits_{e \in E}[[0, |X \cap X_e|]]$. We have 
\begin{align*}
    (R_\rho \circ \pi_\rho)(X) & = \max\limits_{ \textbf{b} \in I_\rho^X} \left(\sum\limits_{e \in E}b_e\right) .
\end{align*}
For any $\textbf{b} \in I_\rho^X$, it follows that 
\begin{align*}
    \sum\limits_{e \in A}b_e &\leq \rho(A), \text{ for all $A \subseteq E$, and}\\
    b_e &\leq |X \cap X_e| \text{ for all $e \in E$}.
\end{align*}
This implies 
\begin{align*}
    \sum\limits_{e \in E}b_e \leq \min\limits_{B \subseteq A} \left(\rho(B) + \sum\limits_{e \in A-B}|X \cap X_e|\right) .
\end{align*}
We will show that the following function $\phi$ is monotone and submodular.
\begin{align*}
    \phi:2^E &\to \mathbb{Z}_{\geq 0}\\
    A &\mapsto \min\limits_{B \subseteq A} \left(\rho(B) + \sum\limits_{e \in A-B}|X \cap X_e|\right) 
\end{align*}
For monotonicity, let $A_1, A_2 \in 2^E$ with $A_1 \subseteq A_2$. Let $B' \subseteq A_2$ be such that 
\begin{align*}
    \rho(B') + \sum\limits_{e \in A_2-B'}|X \cap X_e| &= \min\limits_{B \subseteq A_2} \left(\rho(B) + \sum\limits_{e \in A_2-B}|X \cap X_e|\right) \\
    & = \phi(A_2).
\end{align*}
Since $\rho$ is a polymatroid and thus monotone, we have $\rho(B' \cap A_1) \leq \rho(B')$. Since $A_1 \subseteq A_2$, we have 
\begin{align*}
    \sum\limits_{e \in A_1 - (B'\cap A_1)}|X \cap X_e| \leq \sum\limits_{e \in A_2-B'}|X \cap X_e|.
\end{align*}
Thus 
\begin{align*}
    \phi(A_1) \leq \rho(B' \cap A_1) + \sum\limits_{e \in A_1 - (B'\cap A_1)}|X \cap X_e| \leq\phi(A_2),
\end{align*}
as desired. In particular, this ensures $
\phi(A) \leq \phi(E)$ for any $A \subseteq E$. For submodularity, let $f,g \in E-A$. Consider the following square of values of $\phi$:
$$\begin{tikzcd}
\substack{\phi(A \cup \{g\}) =  \\\min\limits_{B \subseteq A \cup \{g\}} \left(\rho(B) + \sum\limits_{e \in (A \cup \{g\})-B}|X \cap X_e|\right) }\arrow[d, no head] \arrow[r, no head] & \substack{\phi(A \cup \{f,g\}) =  \\\min\limits_{B \subseteq A \cup \{f,g\}} \left(\rho(B) + \sum\limits_{e \in (A \cup \{f,g\})-B}|X \cap X_e|\right) }\arrow[d, no head] \\
\substack{\phi(A) =  \\\min\limits_{B \subseteq A} \left(\rho(B) + \sum\limits_{e \in A-B}|X \cap X_e|\right) }                                         & \substack{\phi(A \cup \{f\}) =  \\\min\limits_{B \subseteq A \cup \{f\}} \left(\rho(B) + \sum\limits_{e \in (A \cup \{f\})-B}|X \cap X_e|\right) } \arrow[l, no head] 
\end{tikzcd}$$
Let $B_f \subseteq A \cup \{f\}$ be such that
\begin{align*}
    \rho(B_f) + \sum\limits_{e \in (A \cup \{f\}-B_f)}|X \cap X_e| &= \min\limits_{B \subseteq A \cup \{f\}} \left(\rho(B) + \sum\limits_{e \in (A \cup \{f\})-B}|X \cap X_e|\right) \\
    & = \phi(A \cup \{f\})
\end{align*}
and define $B_g \subseteq A \cup \{g\}$ similarly. Since $\rho$ is submodular, by definition,
\begin{align*}
    \rho(B_f) + \rho(B_g) \geq \rho(B_f \cap B_g) + \rho(B_f \cup B_g).
\end{align*}
We also have
\begin{align*}
    \sum\limits_{e \in A \cup \{f\}}|X \cap X_e| + \sum\limits_{e \in A \cup \{g\}}|X \cap X_e| = \sum\limits_{e \in A}|X \cap X_e| + \sum\limits_{e \in A \cup \{f,g\}}|X \cap X_e|.
\end{align*}
Therefore, since $B_f \cap B_g \subseteq A$ and $B_f \cup B_g \subseteq A \cup \{f,g\}$,
\begin{align*}
    \phi(A \cup \{f\}) + \phi(A \cup \{g\}) &\geq \rho(B_f \cap B_g) + \sum\limits_{e \in A-(B_f \cap B_g)}|X \cap X_e|\\
    &\hspace{8ex}+ \rho(B_f \cup B_g) + \sum\limits_{e \in A -(B_f \cup B_g)}|X \cap X_e|\\
    &\geq \phi(A) + \phi(A \cup \{f,g\}).
\end{align*}
Monotonicity and submodularity of $\phi$ imply tightness of the bound
\begin{align*}
    \sum\limits_{e \in E}b_e \leq \min\limits_{B \subseteq E}\left(\rho(B) + \sum\limits_{E - B} |X \cap X_e|\right)  .
\end{align*}
We have shown that the largest possible coordinate sum of a vector in $I_\rho^X$ equals $\phi(E)$, i.e. $R(X) = r_{\nat }(X)$ for all $X \subseteq X_E$.
\end{proof}

\begin{example}
    Consider again the $3$-polymatroid $\rho$ from Example \ref{ex-{example1}}. We analyze its $2$-compression $\tau = \rho_{\downarrow e}^2$. We have $X_{E-e} = \{f_1,f_2,f_3\}$. Then
    \begin{align*}
        r_{M_\tau^3}(\varnothing) &= 0\\
        r_{M_\tau^3}(\{f_i\}) &= 1\\
        r_{M_\tau^3}(\{f_i, f_j\}) &= r_{M_\tau^3}(\{f_1, f_2, f_3\}) = 2.
    \end{align*}
In particular, note that
\begin{align*}
    \rho^2_{\downarrow e}(\varnothing) &= r_{M_\tau^3}(\varnothing) = 0\\
    \rho^2_{\downarrow e}(\{f\}) &= r_{M_\tau^3}(\{f_i, f_j\}) = r_{M_\tau^3}(\{f_1, f_2, f_3\}) = 2.
\end{align*}    
$$\begin{tikzpicture}
		[ lightdashed/.style={very thin,gray, dashed},
            darkdashed/.style={very thick, black, dashed},
            darksolid/.style = {very thick, black},
			axis/.style={->,black, very thin}]
			
	\draw[axis] (0,0)--(4,0) node[anchor=west]{$e$};
	\draw[axis] (0,0) -- (0,4) node[anchor=west]{$f$};

    \draw[darksolid, fill = white] (0,0)--(3,0)--(3,1)--(2,2)--(0,2)--cycle;

    \draw[lightdashed] (1,0)--(1,4);
    \draw[lightdashed] (2,0)--(2,4);
    \draw[lightdashed] (3,0)--(3,4);
    \draw[lightdashed] (0,1)--(4,1);
    \draw[lightdashed] (0,2)--(4,2);
    \draw[lightdashed] (0,3)--(4,3);

    \node[] at (-0.2,-0.3) {\tiny$\substack{R_\rho(0,0) \\= 0}$}; 
    \node[] at (-0.2,0.8) {$1$}; 
    \node[] at (0.8,-0.2) {$1$}; 
    \node[] at (0.8,0.8) {$2$}; 
    \node[] at (1.8,-0.2) {$2$}; 
    \node[] at (-0.2, 1.8) {$2$}; 
    \node[] at (2.8,-0.2) {$3$}; 
    \node[] at (1.8,0.8) {$3$}; 
    \node[] at (0.8,1.8) {$3$}; 
    \node[] at (-0.2, 2.8) {$2$};
    \node[] at (2.8,0.8) {$4$}; 
    \node[] at (1.8,1.8) {$4$}; 
    \node[] at (0.8,2.8) {$3$}; 
    \node[] at (1.8,2.8) {$4$}; 
    \node[] at (2.8,2.8) {$4$}; 
    \node[] at (2.8,1.8) {$4$}; 

    \draw[] (1.6,-0.4)--(2.4,-0.4)--(2.4,3.1)--(1.6,3.1)--cycle;

    \node[shape = circle, fill = black, scale = 0.5] at (2,0) {};
    \node[] at (2.2,0.2) {$\textbf{x}$}; 


	\draw[axis] (5,0)--(9,0) node[anchor=west]{$e$};
	\draw[axis] (5,0) -- (5,4) node[anchor=west]{$f$};

    \draw[darksolid, fill = white] (5,0)--(8,0)--(8,1)--(7,2)--(5,2)--cycle;

    \draw[lightdashed] (6,0)--(6,4);
    \draw[lightdashed] (7,0)--(7,4);
    \draw[lightdashed] (8,0)--(8,4);
    \draw[lightdashed] (5,1)--(9,1);
    \draw[lightdashed] (5,2)--(9,2);
    \draw[lightdashed] (5,3)--(9,3);

    \draw[] (6.6,-0.4)--(7.4,-0.4)--(7.4,3.1)--(6.6,3.1)--cycle;

    \node[] at (6.8,-0.2) {$0$}; 
    \node[] at (6.8,0.8) {$1$}; 
    \node[] at (6.8,1.8) {$2$}; 
    \node[] at (6.8,2.8) {$2$}; 
 \end{tikzpicture}$$
 Left: We take a `slice' of $L_\rho^3$ at $e = l = 2$. Right: We subtract $R_\rho(\textbf{x}) = 2$ from each value in the slice, revealing the values of $r_{M_\tau^3}$ (and therefore $\rho_{\downarrow e}^2$).
\end{example}

Since $\rho_{\backslash \{e\}} = \rho_{\downarrow e}^0$ and $\rho_{/\{e\}} = \rho_{\downarrow e}^{m}$ where $m \geq \rho(\{e\})$, the data of the deletions and contractions of $\rho$ can be easily deduced from the values of $R_\rho$ on  $\{0,k\}^E$.

\begin{example} We continue with the $3$-polymatroid $\rho$ from Example \ref{ex-{example1}}. The values of $R_\rho$ on $\{0,3\}^E$ are boxed:
$$
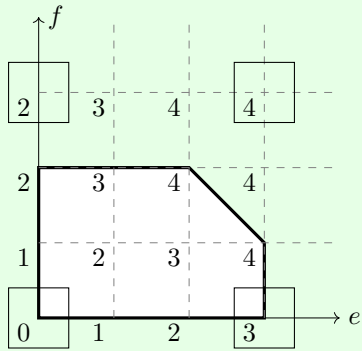
\begin{tikzpicture}
		[ lightdashed/.style={very thin,gray, dashed},
            darkdashed/.style={very thick, black, dashed},
            darksolid/.style = {very thick, black},
			axis/.style={->,black, very thin}]
			
	\draw[axis] (0,0)--(4,0) node[anchor=west]{$e$};
	\draw[axis] (0,0) -- (0,4) node[anchor=west]{$f$};

    \draw[darksolid, fill = white] (0,0)--(3,0)--(3,1)--(2,2)--(0,2)--cycle;

    \draw[lightdashed] (1,0)--(1,4);
    \draw[lightdashed] (2,0)--(2,4);
    \draw[lightdashed] (3,0)--(3,4);
    \draw[lightdashed] (0,1)--(4,1);
    \draw[lightdashed] (0,2)--(4,2);
    \draw[lightdashed] (0,3)--(4,3);

    \node[] at (-0.2,-0.2) {$0$}; 
    \node[] at (-0.2,0.8) {$1$}; 
    \node[] at (0.8,-0.2) {$1$}; 
    \node[] at (0.8,0.8) {$2$}; 
    \node[] at (1.8,-0.2) {$2$}; 
    \node[] at (-0.2, 1.8) {$2$}; 
    \node[] at (2.8,-0.2) {$3$}; 
    \node[] at (1.8,0.8) {$3$}; 
    \node[] at (0.8,1.8) {$3$}; 
    \node[] at (-0.2, 2.8) {$2$};
    \node[] at (2.8,0.8) {$4$}; 
    \node[] at (1.8,1.8) {$4$}; 
    \node[] at (0.8,2.8) {$3$}; 
    \node[] at (1.8,2.8) {$4$}; 
    \node[] at (2.8,2.8) {$4$};  
    \node[] at (2.8,1.8) {$4$}; 

    \draw[] (-0.4,-0.4)--(0.4,-0.4)--(0.4,0.4)--(-0.4,0.4)--cycle;
    \draw[] (2.6,-0.4)--(3.4,-0.4)--(3.4,0.4)--(2.6,0.4)--cycle;
    \draw[] (-0.4,2.6)--(0.4,2.6)--(0.4,3.4)--(-0.4,3.4)--cycle;
    \draw[] (2.6,2.6)--(3.4,2.6)--(3.4,3.4)--(2.6,3.4)--cycle;

    \end{tikzpicture}$$

The values of the minors $\rho'$ of $\rho$ are:
$$
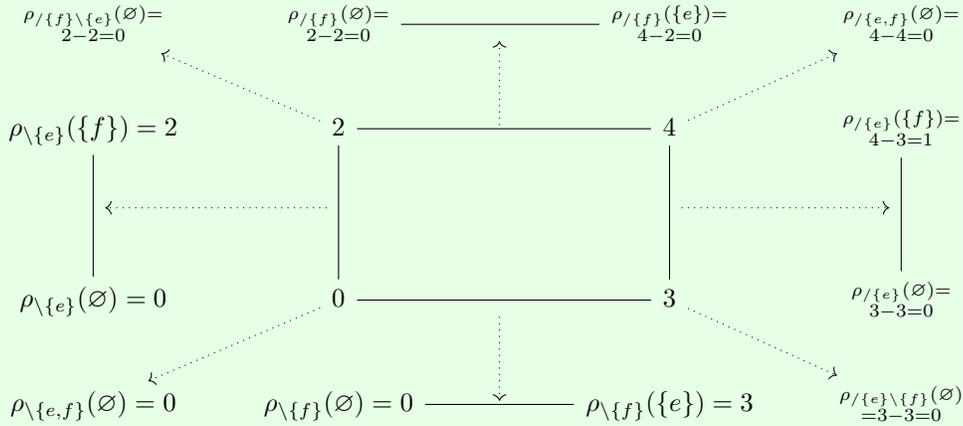
\begin{tikzcd}
\substack{\rho_{/\{f\}\backslash \{e\}}(\varnothing) = \\2-2 = 0} & \substack{\rho_{/\{f\}}(\varnothing) = \\2-2 = 0} \arrow[rr, no head] & {}                   & \substack{\rho_{/\{f\}}(\{e\}) = \\4-2 = 0}      & \substack{\rho_{/\{e,f\}}(\varnothing) = \\4-4 = 0}            \\
\rho_{\backslash \{e\} }(\{f\}) = 2 \arrow[dd, no head]         & 2 \arrow[dd, no head] \arrow[rr, no head] \arrow[lu, dotted]    & {} \arrow[u, dotted] & 4 \arrow[dd, no head] \arrow[ru, dotted] & \substack{\rho_{/\{e\}}(\{f\}) = \\4-3 = 1} \arrow[dd, no head] \\
{}                                                      & {} \arrow[l, dotted]                                            &                      & {} \arrow[r, dotted]                     & {}                                                      \\
\rho_{\backslash \{e\}}(\varnothing) = 0                      & 0 \arrow[rr, no head] \arrow[ld, dotted]                        & {} \arrow[d, dotted] & 3 \arrow[rd, dotted]                     & \substack{\rho_{/\{e\}}(\varnothing) = \\3-3 = 0}             \\
\rho_{\backslash \{e,f\}}(\varnothing) = 0                     & \rho_{\backslash \{f\}}(\varnothing) = 0 \arrow[rr, no head]          & {}                   & \rho_{\backslash \{f\}}(\{e\}) = 3               & \substack{\rho_{/\{e\}\backslash \{f\}}(\varnothing)  \\=3-3 = 0}
\end{tikzcd}$$
\end{example}

\begin{definition}
    Let $(E,\rho)$ be a polymatroid. Fix a minor $\rho'$ of $\rho$, where $\rho' = \rho_{/A_1\backslash A_2}$. Let $H_1$ be the intersection of all hyperplanes $e = \rho(\{e\})$ where $e \in A_1$ and let $H_2$ be the intersection in $\mathbb{R}^E$ of all hyperplanes $e = 0$ where $e \in A_2$. Define the following
    \begin{align*}
        \widetilde{I}_{\rho'} &:= H_1 \cap H_2 \cap I_\rho\\
        \widetilde{B}_{\rho'} &:= \left\{\textbf{x} \in \widetilde{I}_{\rho'}: \text{ for all } \textbf{y} \in \widetilde{I}_{\rho'}, \sum\limits_{e \in E}x_e \geq \sum\limits_{e \in E}y_e\right\} \\
        F_{\rho'} & := H_1 \cap H_2 \cap [[0,k]]^E.
    \end{align*}
\end{definition}

The polytopes $\widetilde{I}_{\rho'}$ and $I_{\rho'}$ are equivalent up to translation, as are $\widetilde{B}_{\rho'}$ and $B_{\rho'}$. The smallest face of the $k$-cube $[[0,k]]^E$ containing $\widetilde{I}_{\rho'}$ is $F_{\rho'}$.

\begin{example}\label{ex-{3D example}}
Consider the $3$-polymatroid $\rho$ on $E = \{e, f, g\}$ such that the rank of each singleton is $3$, the rank of each doubleton is $5$, and the total rank is $5$. The independent set polytope of $\rho$ is depicted below with its visible faces labeled. The hidden faces are:  $\widetilde{I}_{\rho_{\backslash \{e\}}}$ (back facet), $\widetilde{I}_{\rho_{\backslash \{f\}}}$ (left facet),  $\widetilde{I}_{\rho_{\backslash \{g\}}}$ (bottom facet), $\widetilde{I}_{\rho_{\backslash \{f,g\}}}$ (edge along $e$-axis), $\widetilde{I}_{\rho_{\backslash \{e,g\}}}$ (edge along $f$-axis), $\widetilde{I}_{\rho_{\backslash \{e,f\}}}$ (edge along $g$-axis), and $\widetilde{I}_{\rho_{\backslash \{e,f,g\}}}$ (origin). Finally, $\widetilde{I}_\rho = I_\rho$.
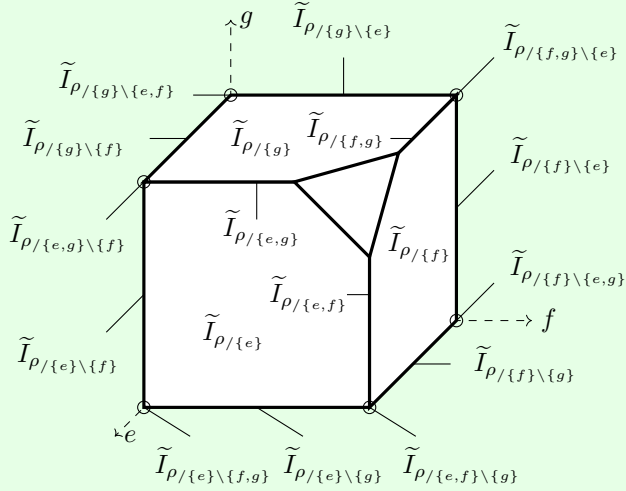
$$\begin{tikzpicture}
		[ lightdashed/.style={very thin,gray, dashed},
            darkdashed/.style={very thick, black, dashed},
            darksolid/.style = {very thick, black},
			axis/.style={->,black, very thin, dashed}]
			
	\draw[axis] (0,0,0) -- (4,0,0) node[anchor=west]{$f$};
	\draw[axis] (0,0,0) -- (0,4,0) node[anchor=west]{$g$};
	\draw[axis] (0,0,0) -- (0,0,4) node[anchor=west]{$e$};

    \draw[darksolid, fill=white!30] (0,0,3) -- (0,3,3) -- (3,3,3) -- (3,0,3) -- cycle;

     \draw[darksolid, fill=white!30] (3,0,0) -- (3,0,3) -- (3,3,3) -- (3,3,0) -- cycle;
    
     \draw[darksolid, fill=white!30] (0,3,0) -- (3,3,0) -- (3,3,3) -- (0,3,3) -- cycle;

     \draw[darksolid, fill=white!30] (2,3,3)--(3,2,3)--(3,3,2)--cycle;

    \node[] at (1.2,1,3) {$\widetilde{I}_{\rho_{/\{e\}}}$};    
    \node[] at (1,3,1.5) {$\widetilde{I}_{\rho_{/\{g\}}}$};
    \node[] at (1,1.8,1.5) {$\widetilde{I}_{\rho_{/\{e,g\}}}$}; 
    \draw[] (1.5,3,3)--(1.5,2.5,3);
    \node[] at (3.1,1.6,1.5) {$\widetilde{I}_{\rho_{/\{f\}}}$};    
    \draw[] (0,3,1.5)--(-0.5,3,1.5);
    \node[] at (-1.5,3,1.5) {$\widetilde{I}_{\rho_{/\{g\}\backslash \{f\}}}$};  
    \draw[] (1.5,3,0)--(1.5,3.5,0);
    \node[] at (1.5,4,0) {$\widetilde{I}_{\rho_{/\{g\}\backslash \{e\}}}$};  
    \draw[] (3,1.5,0)--(3.5,2,0);
    \node[] at (4.4,2.2,0) {$\widetilde{I}_{\rho_{/\{f\}\backslash \{e\}}}$};
    \draw[] (3,0,1.5)--(3.5,0,1.5);
    \node[] at (4.5,0,1.5) {$\widetilde{I}_{\rho_{/\{f\}\backslash \{g\}}}$};
    \draw[] (3, 3, 1.5)--(2.7,3,1.5);
    \node[] at (2, 3,1.2){$\widetilde{I}_{\rho_{/\{f,g\}}}$};

    \draw[] (3,1.5,3)--(2.7, 1.5, 3);
    \node[] at (2.2, 1.5, 3) {$\widetilde{I}_{\rho_{/\{e,f\}}}$};
    
    \draw[] (3,3,0)--(3.5,3.5,0);
    \node[] at (4.4,3.7,0) {$\widetilde{I}_{\rho_{/\{f,g\}\backslash \{e\}}}$};
    \node[shape = circle, draw = black, scale = 0.5] at (3,3,0) {}; 

    \draw[] (3,0,0)--(3.5,0.5,0);
    \node[] at (4.5,0.7,0) {$\widetilde{I}_{\rho_{/\{f\}\backslash \{e,g\}}}$};
    \node[shape = circle, draw = black, scale = 0.5] at (3,0,0) {}; 
    
    \draw[] (3,0,3)--(4,0,4);
    \node[] at (5,0,5) {$\widetilde{I}_{\rho_{/\{e,f\}\backslash \{g\}}}$};
    \node[shape = circle, draw = black, scale = 0.5] at (3,0,3) {}; 

    \draw[] (1.5,0,3)--(2.5,0,4);
    \node[] at (3.3,0,5) {$\widetilde{I}_{\rho_{/\{e\}\backslash \{g\}}}$};

    \draw[] (0,0,3)--(1,0,4);
    \node[] at (1.7,0,5) {$\widetilde{I}_{\rho_{/\{e\}\backslash \{f,g\}}}$};
    \node[shape = circle, draw = black, scale = 0.5] at (0,0,3) {}; 

    \draw[] (0,3,3)--(-0.5,2.5,3);
    \node[] at (-1,2.3,3) {$\widetilde{I}_{\rho_{/\{e,g\}\backslash \{f\}}}$};
    \node[shape = circle, draw = black, scale = 0.5] at (0,3,3) {}; 

    \draw[] (0,1.5,3)--(-0.5,1,3);
    \node[] at (-1,0.7,3) {$\widetilde{I}_{\rho_{/\{e\}\backslash \{f\}}}$};  

    \draw[] (0,3,0)--(-0.5,3,0);
    \node[] at (-1.5,3.2,0) {$\widetilde{I}_{\rho_{/\{g\}\backslash \{e,f\}}}$};
    \node[shape = circle, draw = black, scale = 0.5] at (0,3,0) {}; 
 
    \end{tikzpicture}$$
    
\end{example}

\section{Essential boundedness for polymatroids}
\label{s-{Essential boundedness for polymatroids}}

The singleton and doubleton minors of a polymatroid encode a surprising amount of information about its structural complexity.

\begin{definition}\label{def-{maximally-separated}}
    A matroid is \textit{maximally-separated} if it can be expressed as a direct sum of loops and coloops.
\end{definition}

\begin{definition}
    Let $E$ be a finite set and let $H_1 \subseteq \mathbb{R}^E$ be the unit hypercube given by the Cartesian product $[[0,1]]^{E}$. Consider any hypercube $H$ given by
    \begin{align*}
        H = aH_1 + \textbf{v}
    \end{align*} 
    where $a \geq 0$ and $\textbf{v} \in \mathbb{R}^E_{\geq 0}$. Let $\textbf{u}_e$ denote the unit vector containing $1$ in the coordinate corresponding to $e \in E$ and $0$ everywhere else. For $n < \frac{a}{2}$, we define the \textit{set of $n$-corners of $H$} to be
    \begin{align*}
        Cor_n(H): = \left\{H': H' = (nH_1 + \textbf{v}) + \sum\limits_{e \in A \subseteq E}(a-2n)\textbf{u}_e\right\}
    \end{align*}
    and the \textit{set of corners of $H$} to be
    \begin{align*}
        Cor(H):= \bigcup\limits_{0 \leq n \leq \frac{a}{2}} Cor_n(H).
    \end{align*}
    Intuitively, $Cor_n(H)$ is the set of hypercubic `corners' of $H$ with edge length $n$. The set $Cor(H)$ is the set of hypercubic `corners' of $H$ for all $n$ such that no two elements of $Cor_n(H)$ have a nontrivial intersection in $\mathbb{R}^E$. For simplicity, we sometimes specify an $n$-corner by directly referring to it as a Cartesian product. 
\end{definition}

\begin{definition}\label{def-{n-corner decomposition}}
Let $(E,\rho)$ be a $k$-polymatroid. There must exist some $n \in [0,k]$ such that for some $n$-polymatroid $(E,\tau)$ and maximally-separated matroid $(E,r)$, we can decompose $\rho$ as
\begin{align}\label{eq-{corner decomposition}}
    \rho= \tau+ (k-n)r.
\end{align}
We refer to the right-hand side of Equation \ref{eq-{corner decomposition}} as an \textit{$n$-corner decomposition of $\rho$}. 

From now on, assume $k \geq 2n+1$ unless otherwise stated. This bound on $n$ ensures that the $n$-corner decomposition of $\rho$ is unique: if $\rho(\{e\}) > n$, then $r(\{e\}) = 1$, otherwise $r(\{e\}) = 0$. Let $\cal{N}$ be the set of $n$ for which $\rho$ has an $n$-corner decomposition. Let $m =\min\left(\cal{N}\right)$. We say $\rho$ is \textit{essentially $m$-bounded} and we will refer to its $m$-corner decomposition as the \textit{canonical decomposition of $\rho$}. 
\end{definition}

\begin{example}\label{ex-{confinement to a corner}}
    The $3$-polymatroid $(E,\rho)$ from Example \ref{ex-{3D example}} is essentially $1$-bounded, with canonical decomposition
    \begin{displaymath}
        \rho = U_{2,3} + (3-1)(\rho|_{\{e\}} \oplus \rho|_{\{f\}} \oplus \rho|_{\{g\}}).
    \end{displaymath}
We observe that for each proper minor $\rho'$ of $\rho$, $\widetilde{B}_{\rho'}$ is confined to a $1$-corner of $F_{\rho'}$. This `forces' $B_\rho$ into a $1$-corner of $F_\rho$. Below we display $\widetilde{I}_{\rho'}$ where $\rho'$ is the contraction on $e$, the deletion on $e$, the restriction to $e$, and the contraction to $e$, respectively. The operations on $f$ and $g$ are analogous by symmetry.
$$
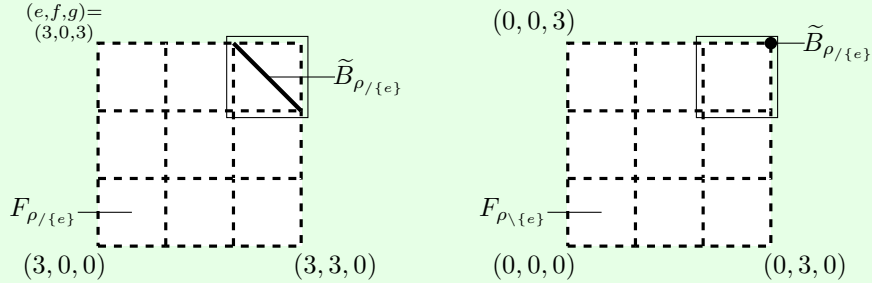
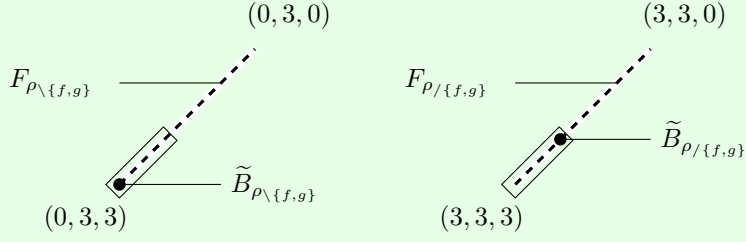
\begin{tikzpicture}
		[ scale = 0.9, lightdashed/.style={very thin,gray, dashed},
            darkdashed/.style={very thick, black, dashed},
            darksolid/.style = {ultra thick, black},
            whiteedge/.style = {line width = 1mm, white},
			axis/.style={->,black, very thin}]

    \draw[darkdashed, fill = white!30] (0,0)--(0,3)--(3,3)--(3,0)--cycle;
    \draw[] (1.9,1.9)--(3.1,1.9)--(3.1,3.1)--(1.9,3.1)--cycle;
    \draw[darkdashed, fill = white!30] (1,0)--(1,3);
    \draw[darkdashed, fill = white!30] (2,0)--(2,3);
    \draw[darkdashed, fill = white!30] (0,1)--(3,1);
    \draw[darkdashed, fill = white!30] (0,2)--(3,2);

    \draw[darksolid] (3,2) -- (2,3);
    \node[] at (-0.8,0.5) {$F_{\rho_{/\{e\}}}$}; 
    \draw[] (0.5,0.5)--(-0.3,0.5);
    \node[] at (-0.5,3.3) {$\substack{(e,f,g) = \\(3,0,3)}$}; 
    \node[] at (-0.5, -0.3) {$(3,0,0)$}; 
    \node[] at (3.5,-0.3) {$(3,3,0)$}; 

    \draw[](2.5,2.5)--(3.5,2.5);
    \node[] at (4,2.5) {$\widetilde{B}_{\rho_{/\{e\}}}$}; 
    \end{tikzpicture}
    \hspace{5ex}
\begin{tikzpicture}
		[ scale = 0.9, lightdashed/.style={very thin,gray, dashed},
            darkdashed/.style={very thick, black, dashed},
            darksolid/.style = {ultra thick, black},
            whiteedge/.style = {line width = 1mm, white},
			axis/.style={->,black, very thin}]

    \draw[darkdashed, fill = white!30] (0,0)--(0,3)--(3,3)--(3,0)--cycle;
    \draw[] (1.9,1.9)--(3.1,1.9)--(3.1,3.1)--(1.9,3.1)--cycle;
    \draw[darkdashed, fill = white!30] (1,0)--(1,3);
    \draw[darkdashed, fill = white!30] (2,0)--(2,3);
    \draw[darkdashed, fill = white!30] (0,1)--(3,1);
    \draw[darkdashed, fill = white!30] (0,2)--(3,2);

    \node[] at (-0.8,0.5) {$F_{\rho_{\backslash \{e\}}}$}; 
    \draw[] (0.5,0.5)--(-0.3,0.5);
    \node[] at (-0.5,3.3) {$(0,0,3)$}; 
    \node[] at (-0.5, -0.3) {$(0,0,0)$}; 
    \node[] at (3.5,-0.3) {$(0,3,0)$}; 
    \node[shape = circle, fill = black, scale = 0.5] at (3,3) {};

    \draw[](3,3)--(3.5,3);
    \node[] at (4,3) {$\widetilde{B}_{\rho_{/\{e\}}}$}; 

 \end{tikzpicture}$$
    
    $$\begin{tikzpicture}
    [ scale = 0.9, lightdashed/.style={very thin,gray, dashed},
            darkdashed/.style={very thick, black, dashed},
            darksolid/.style = {ultra thick, black},
            whiteedge/.style = {line width = 1mm, white},
			axis/.style={->,black, very thin}]

    \draw[whiteedge] (0,0)--(2,2);
    \draw[darkdashed] (0,0)--(2,2);

    \draw[] (-0.2,0)--(0,-0.2)--(0.85,0.65)--(0.65,0.85)--cycle;

    \node[shape = circle, fill = black, scale = 0.5] at (0,0) {};

    \draw[] (1.5,1.5)--(0,1.5);
    \node[] at (-1,1.5) {$F_{\rho_{\backslash \{f,g\}}}$}; 
    \draw[] (0,0)--(1.5,0);
    \node[] at (2.3,0) {$\widetilde{B}_{\rho_{\backslash \{f,g\}}}$}; 
    \node[] at (-0.5,-0.5) {$(0,3,3)$}; 
    \node[] at (2.5,2.5) {$(0,3,0)$}; 
    \end{tikzpicture}
    \hspace{5ex}
        \begin{tikzpicture}
    [ scale = 0.9, lightdashed/.style={very thin,gray, dashed},
            darkdashed/.style={very thick, black, dashed},
            darksolid/.style = {ultra thick, black},
            whiteedge/.style = {line width = 1mm, white},
			axis/.style={->,black, very thin}]

    \draw[whiteedge] (0,0)--(2,2);
    \draw[darkdashed] (0,0)--(2,2);

    \draw[] (-0.2,0)--(0,-0.2)--(0.85,0.65)--(0.65,0.85)--cycle;

    \node[shape = circle, fill = black, scale = 0.5] at (2/3,2/3) {};

    \draw[] (1.5,1.5)--(0,1.5);
    \node[] at (-1,1.5) {$F_{\rho_{/ \{f,g\}}}$}; 
    \draw[] (2/3,2/3)--(2,2/3);
    \node[] at (2.8,2/3) {$\widetilde{B}_{\rho_{/\{f,g\}}}$}; 
    \node[] at (-0.5,-0.5) {$(3,3,3)$}; 
    \node[] at (2.5,2.5) {$(3,3,0)$}; 
    \end{tikzpicture}$$
    \end{example}

In fact, it suffices to investigate only the singleton and doubleton minors of a polymatroid $\rho$ to determine the exact corner of $F_\rho$ to which $B_\rho$ is confined. The next proposition makes this idea precise.

\begin{proposition}\label{prop-{singleton and doubleton minors}}
    Let $(E,\rho)$ be a $k$-polymatroid and let $m \in \mathbb{Z}_{\geq 0}$ such that $k \geq 3m+1$. If each singleton and doubleton minor $\rho'$ of $\rho$ has an $m$-corner decomposition, then $\rho$ has an $m$-corner decomposition.
\end{proposition}

\begin{proof}
We induct on $|E|$. The cases $|E| = 1, 2$ are given by the hypotheses of the theorem. Now fix $|E| \geq 3$ and assume the proposition holds for polymatroids on smaller ground sets. Fix $e \in E$. By the inductive hypothesis, we can write
\begin{align*}
    \rho_{\backslash \{e\}}  &= \taudel + (k-m)\rdel\\
    \rho_{/\{e\}} & = \taucont + (k-m)\rcont\\
    \rho|_{\{e\}} & = \taures + (k-m)\rres.
\end{align*}
For $\varphi = \tau, r$, define $\varphi: 2^E \to \mathbb{Z}_{\geq 0}$ as follows,
\begin{align*}
    \varphi(A) := \begin{cases}
        \varphi^{\text{del}}(A) & \text{if $e \notin A$,}\\
        \varphi^{\text{res}}(\{e\}) + \varphi^{\text{cont}}(A\backslash \{e\}) & \text{if $e \in A$.}
    \end{cases}
\end{align*}
It is easy to check that
\begin{align*}
    \rho = \tau + (k-m)r.
\end{align*}
\begin{example} Consider the $3$-polymatroid $(E,\rho)$ from Example \ref{ex-{3D example}}. Recall that its canonical decomposition is $\rho = \tau + (3-1)r$ where
\begin{align*}
    \tau = U_{2,3} \text{ and } r = \rho|_{\{e\}} \oplus \rho|_{\{f\}} \oplus \rho|_{\{g\}}.
\end{align*} 
We highlight two subsets of $I_\rho$ that are equivalent to $I_{(3-1)r}$ and $I_\tau$ up to translation.
$$\begin{tikzpicture}
		[scale = 0.9, lightdashed/.style={very thin,gray, dashed},
            darkdashed/.style={very thick, black, dashed},
            darksolid/.style = {very thick, black},
			axis/.style={->,black, very thin}]
   \begin{scope}[rotate around x = 0]
			
	\draw[axis] (0,0,0) -- (4,0,0) node[anchor=west]{$f$};
	\draw[axis] (0,0,0) -- (0,4,0) node[anchor=west]{$g$};
	\draw[axis] (0,0,0) -- (0,0,4) node[anchor=west]{$e$};
    \draw[lightdashed] (0,0,0) -- (0,3,0) -- (3,3,0) -- (3,0,0) -- cycle;
    \draw[lightdashed] (1,3,0)--(1,0,0);
    \draw[lightdashed] (2,3,0)--(2,0,0);
    \draw[lightdashed] (0,1,0)--(3,1,0);
    \draw[lightdashed] (0,2,0)--(3,2,0);

    \draw[lightdashed] (0,0,0) -- (0,0,3) -- (0,3,3) -- (0,3,0) -- cycle;
    \draw[lightdashed] (0,1,0) -- (0,1,3);
    \draw[lightdashed] (0,2,0) -- (0,2,3);
    \draw[lightdashed] (0,0,1) -- (0,3,1);
    \draw[lightdashed] (0,0,2) -- (0,3,2);

    \draw[lightdashed] (0,0,0) -- (3,0,0) -- (3,0,3) -- (0,0,3) -- cycle;
    \draw[lightdashed] (1,0,0) -- (1,0,3);
    \draw[lightdashed] (2,0,0) -- (2,0,3);
    \draw[lightdashed] (0,0,1) -- (3,0,1);
    \draw[lightdashed] (0,0,2) -- (3,0,2);

    \draw[darksolid, fill=white!30] (0,2,2) -- (2,2,2) -- (2,0,2) -- (0,0,2) -- (0,2,2);
    \draw[darksolid, fill=white!30] (0,2,0) -- (0,2,2) -- (2,2,2) -- (2,2,0) -- (0,2,0);
    \draw[darksolid, fill=white!30] (2,2,0) -- (2,0,0) -- (2,0,2) -- (2,2,2) -- (2,2,0);
    \node[] at (.5, 0.5, 1) {$I_{(3-1)r}$}; 
    \draw[darksolid, fill=white!30] (3,3,2) -- (2,3,2) -- (2,3,3);
    \draw[darksolid, fill=white!30] (2,3,3) -- (2,2,3) -- (3,2,3);
    \draw[darksolid, fill=white!30] (3,3,2) -- (3,2,2) -- (3,2,3);
    \draw[darksolid, fill=white!30] (2,3,3) -- (3,2,3) -- (3,3,2) -- (2,3,3);
    \node[] at (2.7,2.7,3) {$I_\tau$}; 

    \draw[lightdashed] (1,1,0) -- (1,1,3);
    \draw[lightdashed] (2,1,0) -- (2,1,3);
    \draw[lightdashed] (0,1,1) -- (3,1,1);
    \draw[lightdashed] (0,1,2) -- (3,1,2);
    \draw[lightdashed] (1,2,0) -- (1,2,3);
    \draw[lightdashed] (2,2,0) -- (2,2,3);
    \draw[lightdashed] (0,2,1) -- (3,2,1);
    \draw[lightdashed] (0,2,2) -- (3,2,2);

    \draw[lightdashed] (3,0,0) -- (3,0,3) -- (3,3,3) -- (3,3,0) -- cycle;
    \draw[lightdashed] (3,1,0) -- (3,1,3);
    \draw[lightdashed] (3,2,0) -- (3,2,3);
    \draw[lightdashed] (3,0,1) -- (3,3,1);
    \draw[lightdashed] (3,0,2) -- (3,3,2);
    
    \draw[lightdashed] (0,3,0) -- (3,3,0) -- (3,3,3) -- (0,3,3) -- cycle;
    \draw[lightdashed] (1,3,0) -- (1,3,3);
    \draw[lightdashed] (0,3,1) -- (3,3,1);
    \draw[lightdashed] (2,3,0) -- (2,3,3);
    \draw[lightdashed] (0,3,2) -- (3,3,2);

    \draw[lightdashed] (0,0,3) -- (0,3,3) -- (3,3,3) -- (3,0,3) -- (0,0,3);
    \draw[lightdashed] (0,1,3) -- (3,1,3);
    \draw[lightdashed] (1,0,3) -- (1,3,3);
    \draw[lightdashed] (0,2,3) -- (3,2,3);
    \draw[lightdashed] (2,0,3) -- (2,3,3);
    \end{scope}
\end{tikzpicture}$$
\end{example}

\begin{lemma}
        The function $r$ is a maximally-separated matroid.
\end{lemma}
\begin{proof}
        
    Since $\rdel$, $\rcont$, and $\rres$ are maximally-separated matroids, it suffices to show $\rdel(\{f\}) = \rcont(\{f\})$ for all $f \neq e$. Consider this square of subsets:
    $$\begin{tikzcd}
    \{e\} \arrow[d, no head] \arrow[r, no head] & \{e,f\} \arrow[d, no head] \\
    \varnothing\arrow[r, no head]                             & \{f\}
    \end{tikzcd}$$

     The values of $\rho$ on this square are:
\begin{align*}
    \rho(\{e\}) & = \taures(\{e\}) + (k-m)\rres(\{e\})\\
    \rho(\{e, f\}) & = \taures(\{e\}) + \taucont(\{f\}) + (k-m)\rres(\{e\}) + (k-m)\rcont(\{f\})\\
    \rho(\varnothing) & = 0\\
    \rho(\{f\}) & = \taudel(\{f\}) + (k-m)\rdel(\{f\})
\end{align*}

Now consider $\rho' = \rho|_{\{e,f\}}$ which is a doubleton minor of $\rho$. Applying Proposition \ref{prop-{singleton and doubleton minors}} to $\rho'$ gives
\begin{align*}
    \rho' = \tau' + (k-m)r'
\end{align*}
where $\tau'$ is an $m$-polymatroid and $r'$ is a maximally-separated matroid. The values of $\rho'$ on the square are
\begin{align*}
    \rho'(\{e\}) & = \tau'(\{e\}) + (k-m)r'(\{e\})\\
    \rho'(\{e, f\}) & = \tau'(\{e, f\}) + (k-m)r'(\{e\}) + (k-m)r'(\{f\})\\
    \rho'(\varnothing) & = 0\\
    \rho'(\{f\}) & = \tau'(\{f\}) + (k-m)r'(\{f\}).
\end{align*}
Since $\rho$ and $\rho'$ take the same values on $2^{ef}$, we have
\begin{align*}
    \rho(\{e, f\})-\rho(\{e\}) &= \rho'(\{e, f\})-\rho'(\{e\})\\
    \rho(\{f\})-\rho(\varnothing) & = \rho'(\{f\}) - \rho'(\varnothing)
\end{align*}
which we can expand as
\begin{align*}
    \taucont(\{f\}) +(k-m)\rcont(\{f\}) &= \tau'(\{e, f\})-\tau'(\{e\}) + (k-m)r'(\{f\})\\
    \taudel(\{f\}) + (k-m)\rdel(\{f\}) & =\tau'(\{f\}) + (k-m)r'(\{f\}) 
\end{align*}
Subtracting the second equation from the first and then rearranging, we have
\begin{align*}
    (k-m)(\rcont(\{f\})-\rdel(\{f\})) & = \tau'(\{e, f\}) - \tau'(\{e\}) - \tau(\{f\}) \\
    &\hspace{20ex}+ \taudel(\{f\}) - \taucont(\{f\}).
\end{align*}
Since $\taucont, \taudel$, and $\tau'$ are $m$-polymatroids, we have
\begin{align*}
     -m \leq \taudel(\{f\}) - \taucont(\{f\})\leq m\\
     -m \leq \tau'(\{e, f\}) - \tau'(\{e\}) - \tau(\{f\}) \leq 0
\end{align*}
which implies
\begin{align*}
    \tau'(\{e, f\}) - \tau'(\{e\}) - \tau(\{f\}) + \taudel(\{f\}) - \taucont(\{f\}) \equiv x \pmod{k-m}
\end{align*}
where $x \in [0,2m]$. Since
\begin{align*}
    (k-m)(\rcont(\{f\})-\rdel(\{f\}))\equiv 0 \pmod{k-m}
\end{align*}
and $k-m \geq 2m+1$, it must be that $\rcont(\{f\}) = \rdel(\{f\})$, as desired.
\end{proof}

\begin{lemma}
        The function $\tau$ is an $m$-polymatroid.
\end{lemma}

\begin{proof}
    By definition, $\tau$ takes values in $\mathbb{Z}_{\geq 0}$. Using submodularity of $\rho$, it is not difficult to check that $\tau$ is submodular. We have $\tau(\{e\}) = \taures(\{e\}) \leq m$ and $\tau(\{g\}) = \taudel(\{g\}) \leq m$ for any $g \neq e$. 
    
    It remains to show $\tau$ is monotonic, i.e. $\tau(A \cup \{g\}) \geq \tau(A)$ for all $A \subseteq E$ and $g \in E$. If $g \neq e$, then this follows from the monotonicity of $\taudel$ and $\taucont$. Now assume $g = e$. We induct on $|A|$. Since $\tau(\{e\}) \in [0,m]$, the statement follows easily when $A$ is empty. Now fix a nonempty $A$ and assume the lemma holds for smaller subsets. Choose $f \in A$ and write $A = B \sqcup \{f\}$. Consider the square of subsets
    $$\begin{tikzcd}
B \sqcup \{e\} \arrow[d, no head] \arrow[r, no head] & B \sqcup \{e,f\} \arrow[d, no head] \\
B \arrow[r, no head]                             & B \sqcup \{f\}                    
\end{tikzcd}$$
The values of $\rho$ on this square are
\begin{align*}
    \rho(B \sqcup \{e\}) & = \taures(\{e\})  + \taucont(B) + (k-m)\rres(\{e\}) + (k-m)\rcont(B)\\
    \rho(B \sqcup \{e,f\}) & = \taures(\{e\}) + \taucont(B \sqcup \{f\}) + (k-m)\rres(\{e\}) \\
    &\hspace{30ex}+ (k-m)\rcont(B \sqcup \{f\})\\
    \rho(B) & = \taudel(B) + (k-m)\rdel(B) \\
    \rho(B \sqcup \{f\}) & = \taudel(B\sqcup \{f\}) + (k-m)\rdel(B \sqcup \{f\})
\end{align*}
Let $\rho' = (\rho_{/(B \sqcup \{f\})})|_{\{e\}}$. Applying the inductive hypothesis from the proof of Proposition \ref{prop-{singleton and doubleton minors}} to $\rho'$, if $r(\{e\}) = 0$ then $\rho'(\{e\}) \in [0,m]$, and if $r(\{e\}) = 1$ then $\rho'(\{e\}) \in [k-m,k]$. By definition of contraction, the upper-right value minus the lower-right value in the square is equal to $\rho'(\{e\})-\rho'(\varnothing)$. This gives
   \begin{align*}
       \tau(B \sqcup \{e,f\}) - \tau(B \sqcup \{f\}) \equiv x \pmod{k-m}.
   \end{align*}
   where $x \in [0,m]$. 
    In particular, this implies that one of the following must be true:
    \begin{align*}
        \tau(B \sqcup \{e,f\}) &- \tau(B \sqcup \{f\}) \geq 0, \text{ or}\\
        \tau(B \sqcup \{e,f\}) &- \tau(B \sqcup \{f\}) \leq 2m-k.
    \end{align*}
    Using the definition of $\tau$,
    \begin{align}\label{eq-{eq1}}
        \tau(B \sqcup \{e,f\}) - \tau(B \sqcup \{e\})  & = \taucont(B \sqcup \{f\}) - \taucont(B),\\
        \label{eq-{eq2}}\tau(B \sqcup \{f\}) - \tau(B) & = \taudel(B \sqcup \{f\}) - \taudel(B).
    \end{align}
    Since $\taucont$ and $\taudel$ are $m$-polymatroids, we have
    \begin{align}\label{eq-{ineq1}}
        0 \leq \taucont(B \sqcup \{f\}) - \taucont(B) \leq m\\
        \label{eq-{ineq2}} 0 \leq \taudel(B \sqcup \{f\}) - \taudel(B) \leq m.
    \end{align}
    Furthermore, $\tau(B \sqcup \{e\}) \geq \tau(B)$ by the inductive hypothesis. We subtract equation \ref{eq-{eq1}} from equation \ref{eq-{eq2}}, rearrange, and apply the inequalities \ref{eq-{ineq1}} and \ref{eq-{ineq2}}, arriving at
    \begin{align*}
        \tau(B \sqcup \{e,f\}) - \tau(B \sqcup \{f\})  & \geq -m.
    \end{align*}
    Since $k-m \geq 2m+1$, we conclude 
    \begin{align*}
        \tau(B \sqcup \{e,f\}) &- \tau(B \sqcup \{f\}) \geq 0,
    \end{align*}
    that is, $\tau(A \sqcup \{e\}) \geq \tau(A)$, as desired.
\end{proof}
This concludes the proof of Proposition \ref{prop-{singleton and doubleton minors}}.
\end{proof}

The intuition behind this construction is that $\rho$ can be thought of as the deletion $\rho_{\backslash \{e\}}$, the contraction $\rho_{/\{e\}}$, and some ``gluing data'' between them.

\begin{proposition}
\label{prop-{l-compression is contraction or deletion}}
    Let $(E,\rho)$ be an essentially $m$-bounded $k$-polymatroid with canonical decomposition
    \begin{align*}
        \rho = \tau + (k-m)r.
    \end{align*}
    Then for any $e \in E$ and $m \leq l\leq k-m$, $\rho_{\downarrow e}^l$ is equal to $\rho_{/\{e\}}$ or $\rho_{\backslash \{e\}}$.
\end{proposition}
\begin{proof}
    If $l \geq \rho(\{e\})$, then $\rho_{\downarrow e}^l = \rho_{/\{e\}}$. Now assume $l \leq \rho(\{e\})-m$. Consider any $A\subseteq (E-\{e\})$. Using equation (41), and that $\tau$ is an $m$-polymatroid and $r$ is a maximally-separated matroid, we have
\begin{align*}
    \rho(A) + \rho(\{e\}) - \rho(A\cup \{e\}) & = m\big(\tau(A) + \tau(\{e\}) - \tau(A \cup \{e\})\big)\\
    &\leq m.
\end{align*}
Along with submodularity of $\rho$,
\begin{align*}
    0 \leq \rho(A) + \rho(\{e\}) - \rho(A\cup \{e\})\leq m.
\end{align*}

If $\rho(A) + \rho(\{e\}) - \rho(A\cup \{e\}) = 0$, then $A$ and $e$ are disjoint under $\rho$, so $\rho_{\downarrow e}^l(A) = \rho(A)$ as compressing $e$ leaves $A$ unaffected. If $1 \leq \rho(A) + \rho(\{e\}) - \rho(A\cup \{e\}) \leq m$, then since $l < \rho(\{e\})-m$, contracting $l$ clones of $e$ also leaves $A$ unaffected. Hence, $\rho_{\downarrow e}^l(A) = \rho(A)$. Since $\rho_{\downarrow e}^l(A) = \rho_{\backslash \{e\}}(A)$ for all $A \subseteq (E-\{e\})$, we conclude $\rho_{\downarrow e}^l = \rho_{\backslash \{e\}}$.
\end{proof}

\section{Applications to excluded minor problems}
\label{s-{Applications to excluded minor problems}}

\begin{definition}\label{def-{minor-closed, excluded minors}}
    A class of polymatroids is \textit{minor-closed} if for any $\rho$ in the class, all minors of $\rho$ are also in the class. We can completely characterize a minor-closed polymatroid class by its set of \textit{excluded minors}: those polymatroids that are not in the class, whose proper minors are all in the class.  
\end{definition}

From now on, $\cal{C}$ denotes a minor-closed matroid class, and $Ex(\cal{C})$ is its set of excluded minors. Let $\class{k}$ be the class of $k$-polymatroids whose $k$-natural matroids are in $\cal{C}$. By Observation \ref{obs-{k-natural observations}}(2), $\class{k}$ is also minor-closed, as minors of $k$-polymatroids are $k$-polymatroids. 

\begin{lemma}[\cite{young}]\label{lem-{compression}}
Let $(E, \rho)$ be a $k$-polymatroid in $Ex\big(\class{k}\big)$. Fix $e \in E$ with $\rho(\{e\}) \geq 2$. Any internal compression $\rho_{\downarrow e}^l$ is in $Ex\big(\class{k}\big)$ if and only if $\rho_{\downarrow e}^l \notin \class{k}$.
\end{lemma}

\begin{definition}[\cite{young}]\label{def-{Gamma}}
    Define $\Gamma\big(\class{k}\big)$ to be the set of all $\rho \in Ex\big(\class{k}\big)$ such that $\rho_{\downarrow e}^l \in \class{k}$ for all $e \in E$ and $l \in [1, \rho(\{e\})-1]$.
\end{definition}

\begin{remark}[\cite{young}]\label{rem-{sequence of compressions}}
    If $\rho \in Ex\big(\class{k}\big)$ and $\rho \notin \Gamma\big(\class{k}\big)$, then by Lemma \ref{lem-{compression}}, some sequence of compressions of elements $e_1,\hdots, e_j \in E$ each of rank $2$ or higher starting from $\rho$ eventually yields some $\rho' \in \Gamma\big(\class{k}\big)$. Furthermore, every polymatroid in this sequence is an excluded minor. 
\end{remark}

\begin{lemma}\label{lem-{Gamma implies small E}}
    Let $\cal{C}$ be the class of matroids with $Ex(\cal{C}) = \{\unif , \dunif \}$. If $\rho \in \Gamma\big(\class{k}\big)$, then $|E| \leq b$. 
\end{lemma}
\begin{proof}
    For any $e \in E$, it must be that deleting the entirety of $X_e$ or contracting the entirety of $X_e$ from $\nat $ eliminates all $(\unif )$- and $(\dunif )$-minors of $\nat $. This is because $\nat \backslash X_e = M_{\rho\backslash \{e\}}^k$ and $\nat /X_e = M_{\rho/\{e\}}^k$, but $M_{\rho\backslash \{e\}}^k$ and $M_{\rho/ e}^k$ are both in $\cal{C}$ since $\rho \in Ex\big(\class{k}\big)$ implies $\rho\backslash \{e\}$ and $\rho/\{e\}$ are in $\class{k}$. Thus, if $\rho(\{e\}) = 1$, then at least one element of $X_e$ is in each $(\unif )$- or $(\dunif )$-minor of $\nat $, and if $\rho(\{e\}) \geq 2$, then to get a $(\unif )$- or $(\dunif )$-minor of $\nat $, we must do exactly one of the following:
\begin{enumerate}
    \item Contract at least one (but not all) of the elements of $X_e$; delete the rest.
    \item Have at least one element of $X_e$ in the $(\unif )$- or $(\dunif )$-minor.
\end{enumerate}

Assume (1) applies to $e$, and let $l = |Y|$, so $l \in [1, \rho(\{e\})-1]$. Then the $k$-natural matroid of $\rho_{\downarrow e}^l$ is isomorphic to $\nat /Y\backslash (X_e-Y)$, which by assumption contains a $(\unif )$- or $(\dunif )$-minor. Therefore, it must be that $\rho\notin\Gamma\big(\class{k}\big)$. If $\rho \in \Gamma\big(\class{k}\big)$, then (1) cannot apply to any $e \in E$. Since (2) can occur for at most $b$ elements of $E$, it must be that $|E| \leq b$.
\end{proof}

Proposition \ref{prop-{l-compression is contraction or deletion}} allows us to systematically discard polymatroids which would otherwise be legitimate candidates for excluded minors for certain classes of polymatroids defined by their $k$-natural matroids. 

\begin{theorem}\label{th-{exclude unif rank 2}}
For $b \geq 4$, let $\cal{C}$ be the class of matroids characterized by
    \begin{align*}
        Ex(\cal{C}) = \{U_{2,b}, U_{b-2,b}\}.
    \end{align*}
    The set $Ex\big(\class{k}\big)$ is finite for $k \geq 2b-4$.
\end{theorem}
\begin{proof} \leavevmode
We will show that every polymatroid in $Ex\big(\class{k}\big)$ has ground set size at most $b$; finiteness immediately follows. Assume towards a contradiction that there exists $\rho \in Ex\big(\class{k}\big)$ such that $|E(\rho)| > b$. 

First, we will show that $\rho$ has a $1$-corner decomposition. By Proposition \ref{prop-{singleton and doubleton minors}}, it suffices to show that each singleton and doubleton minor $\rho'$ of $\rho$ has a $1$-corner decomposition.

Consider $|E(\rho')| = 1$ where $E(\rho') = \{e\}$. Any proper minor of $\rho'$ is empty and therefore trivially has a $1$-corner decomposition. Since $\rho$ is an excluded minor, $\rho'$ cannot be an excluded minor. Therefore, using Theorem 3.0.1 in \cite{young} and the inequality $k \geq 2b-4$, we conclude $\rho'(\{e\}) \in \{0, 1, k-1, k\}$. For each choice of $\rho'(\{e\})$, the $1$-corner decomposition $\rho' = \tau + (k-1)r$ is given by:
    \begin{table}[H]
\centering
\begin{tabular}{ |c|c|c| }
 \hline
$\rho'(\{e\})$& $\tau$&$r$\\
\hline
$0$&$U_{0,1}$&$U_{0,1}$\\
 \hline
 $1$&$\coloop$&$U_{0,1}$\\
 \hline
 $k-1$&$U_{0,1}$&$\coloop$\\
 \hline
 $k$&$\coloop$&$\coloop$\\
 \hline
\end{tabular}
\end{table}
Next, consider $|E(\rho')| = 2$ where $E(\rho') = \{e, f\}$. Any proper minor of $\rho'$ is empty or falls into the singleton case so it suffices to address $\rho'$ itself. We use Theorem 4.0.1 from \cite{young} to determine the allowed possibilities for the first three columns in the table below. For each choice of $\rho'$, the $1$-corner decomposition $\rho' = \tau + (k-1)r$ is given by\footnote{Note: Direct sums are written in the form $\rho|_{\{e\}} \oplus \rho|_{\{f\}}$.}:
    \begin{table}[H]
\centering
\begin{tabular}{ |c|c|c|c|c|c| }
 \hline
$\rho'(\{e\})$ & $\rho'(\{f\})$& $\rho'(\{e, f\})$&$\tau$&$r$\\
\hline
$k-1$ & $k-1$& $2k-2$&$U_{0,1}  \oplus U_{0,1} $&$\coloop \oplus \coloop$\\
 \hline
 $k$& $k$& $2k-1$&$U_{1,2}$&$\coloop \oplus \coloop$\\
 \hline
   $k$& $k$&$2k$&$\coloop \oplus \coloop$&$\coloop \oplus \coloop$\\
  \hline
 $1$ & $k-1$ &$k$&$U_{0,1}  \oplus \coloop$&$\coloop \oplus U_{0,1} $\\
 \hline
 $1$ & $k$&$k$&$U_{1,2}$&$U_{0,1}  \oplus \coloop$\\
\hline
 $1$ & $k$&$k+1$&$\coloop \oplus \coloop$&$U_{0,1}  \oplus \coloop$\\
\hline
$k-1$ & $k$&$2k-1$&$\coloop \oplus U_{0,1} $&$\coloop \oplus \coloop$\\
\hline
\end{tabular}
\end{table}
By Lemma \ref{lem-{Gamma implies small E}}, since $|E(\rho)| > b$, $\rho \notin \Gamma\big(\class{k}\big)$. By Remark \ref{rem-{sequence of compressions}}, some interior compression $\rho_{\downarrow e}^l$ must be in $Ex\big(\class{k}\big)$. However, since $\rho$ has a $1$-corner decomposition, by Proposition \ref{prop-{l-compression is contraction or deletion}}, $\rho_{\downarrow e}^l$ is equal to a proper minor of $\rho$, leading to the desired contradiction.
\end{proof}

To conclude, we investigate what happens when we exclude uniform matroids more generally. That is, let $\cal{C}$ be the class of matroids characterized by
    \begin{align*}
        Ex(\cal{C}) = \{\unif, \dunif\}.
    \end{align*}
Although the situation becomes more complicated when $a \geq 3$ and $b \geq 2a$, Proposition \ref{prop-{singleton and doubleton minors}} does impose a bound on the complexity of $k$-polymatroids $\rho \in Ex\big(\class{k}\big)$ when $k$ is sufficiently large. First, we generalize Propositions 3.0.1 and 4.0.1 in \cite{young}.

\begin{proposition}\label{prop-{general unif excluded singletons}}
    For $b \geq 2a$, let $\cal{C}$ be the class of matroids characterized by
    \begin{align*}
        Ex(\cal{C}) = \{\unif, \dunif\}.
    \end{align*}
    Let $k \geq 2(b-a)$. Consider the $k$-polymatroid $(\{e\}, \rho)$ of rank $m$.
    \begin{enumerate}
        \item If $m \in [0, a-1] \cup [k-a+1, \rho(\{e\})]$, then $\rho \in \class{k}$ and hence is not an excluded minor.
        \item If $m \in [a,k-a]$, then $\rho \in Ex\big(\class{k}\big)$. We denote $\rho$ as $Ex^m$.
    \end{enumerate}
    Hence, there are $k-2a+1$ singleton excluded minors for $\class{k}$.
\end{proposition}

\noindent
We start by proving the following lemma.
\begin{lemma}\label{lem-{matroid nullity}} If a matroid $M$ has a $(\unif )$- or $(\dunif )$-minor, then $\eta(M) \geq a$.\end{lemma}
\renewcommand*{\proofname}{Proof of Lemma \ref{lem-{matroid nullity}}.}
\begin{proof}
    Let $n = |E(M)|$ and $r = r(M)$. To obtain a $(\unif )$-minor of $M$, we must decrease the ground set by $n-b$ and decrease the rank by $r-a$. It must be that $n-b \geq r-a$, so $\eta(M) = n-r \geq b-a$. Similarly, to obtain a $(\dunif )$-minor of $M$, it must be that $\eta(M) = n-r \geq a$. Since $b \geq 2a$, $b-a \geq a$, so in both cases we require $\eta(M) \geq a$. 
\end{proof}
\renewcommand*{\proofname}{Proof of Proposition \ref{prop-{general unif excluded singletons}}.}
\begin{proof}\leavevmode
\begin{enumerate}
    \item If $m \in [0,a-1]$, then the rank of $\rho$ is not large enough for $\nat $ to contain a $(\unif )$- or $(\dunif )$-minor, so $\rho \in \class{k}$. Next, let $m \in [k-a+1, \rho(\{e\})]$. The quantity $\eta(\nat )$ is equal to $k-m$, and $a-1 < a$, so by Lemma \ref{lem-{matroid nullity}}, $\nat $ contains neither $\unif $ nor $\dunif $ as a minor. We conclude that $\rho \in \class{k}$. 
    \item Let $X$ and $Y$ be disjoint subsets of $X_e$ with $m-a$ and $k-m+a-b$ elements respectively. Then $\nat /X\backslash Y \cong \unif $, so $\rho\notin \class{k}$. Any proper minor $\rho'$ of $\rho$ is empty, so $M_{\rho'}^k$ contains neither $\unif $ nor $\dunif $ as a minor. We conclude $\rho \in Ex\big(\class{k}\big)$.\qedhere
\end{enumerate}
\end{proof}

\renewcommand*{\proofname}{Proof.}

We briefly introduce some information that will be used in the proof of Proposition \ref{prop-{general unif excluded doubletons}}.

\begin{lemma}[\cite{young}]\label{lem-{simplification}}
    Let $\mathcal{C}$ be a minor-closed class of matroids whose excluded minors are all simple. If the polymatroid $\rho$ is an excluded minor for $\class{k}$, then $\rho$ cannot have any loops (elements of rank $0$) or nontrivial parallel classes of points (elements of rank $1$).
\end{lemma}

\begin{definition}[\cite{young}]
    Let $(E,\rho)$ be a polymatroid. We will let $S(\rho)$ denote the \textit{simplification} of $\rho$, defined as the polymatroid obtained by deleting all loops of $\rho$, and in each nontrivial parallel class of points of $\rho$, deleting all points except for one.
\end{definition}

If $\rho \neq S(\rho)$, then $S(\rho)$ is a proper minor of $\rho$, so if $S(\rho)$ is not in $\class{k}$, then neither is $\rho$. When analyzing matroid minors of $\nat$ to see if any of them are isomorphic to the matroids in $Ex(\mathcal{C})$, Lemma \ref{lem-{simplification}} allows us to consider minors of $S(\nat)$ instead.

\begin{proposition}\label{prop-{general unif excluded doubletons}}
    For $b \geq 2a$, let $\cal{C}$ be the class of matroids characterized by
    \begin{align*}
        Ex(\cal{C}) = \{\unif , \dunif \}.
    \end{align*}
    Let $k \geq 2(b-a)$. Below, we classify all $k$-polymatroids $\rho = \rho_{(\rho_e,\rho_f)}^m$ on $\{e,f\}$ where $\rho_e = \rho(\{e\})$, $\rho_f = \rho(\{f\})$, and $m = \rho(\{e,f\})$. Let $\rho_e \leq \rho_f$. (A blank cell indicates that there are no additional restrictions on the range of the corresponding value.)
    \label{tab-{doubleton excluded minors for unif}}
    \begin{table}[H]
\centering
\begin{tabular}{ |c|c|c|c|c| }
 \hline
&$\rho_e$ range& $\rho_f$ range&$m$ range & In $\class{k}$ or $Ex\big(\class{k}\big)$?\\
\hline
$1$&$[0,a-1]$&$[0,a-1]$&$[\rho_f, \rho_e+\rho_f]$&$ \class{k}$\\
\hline
$2$&$[0,a-1]$&$[k-a+1, k]$&$[\rho_e + (k-a+1), \rho_e+\rho_f]$&$ \class{k}$\\
\hline
$3$&$[1,a-1]$&$[k-a+1, k]$&$[\rho_f, \rho_e + (k-a)]$&$Ex\big(\class{k}\big)$\\
\hline
$4$&$[k-a+1, k]$&$[k-a+1, k]$&$[\rho_f + (k-a+1), \rho_e + \rho_f]$&$\class{k}$\\
\hline
$5$&$[k-a+1, k]$&$[k-a+1, k]$&$[\rho_f, \rho_f + (k-a)]$&$Ex\big(\class{k}\big)$\\
\hline
$6$&&$[a,k-a]$&& Neither  \\
\hline
$7$&$[a,k-a]$&&& Neither  \\
\hline
\end{tabular}
\end{table}
\end{proposition}
\begin{proof}
    For rows (6) and (7), $\rho|_{\{e\}}$ and $\rho|_{\{f\}}$ respectively are in $Ex\big(\class{k}\big)$, so $\rho$ cannot be in $\class{k}$ nor $Ex\big(\class{k}\big)$. For rows (1) through (5), by Proposition \ref{prop-{general unif excluded singletons}}, no proper minor of $\rho$ contains a $(\unif )$- or $(\dunif )$-minor, so $\rho \in \class{k}$ or $\rho \in Ex\big(\class{k}\big)$.
    \begin{enumerate}
        \item[(1), (4)] Assume by way of contradiction that $\nat $ contains a $(\unif )$- or $(\dunif )$-minor $M$. Let $\nat /X^{\text{cont}}\backslash X^{\text{del}} = M$. Let $X_e^{\text{cont}} = X^{\text{cont}} \cap X_e$ and define $X_e^{\text{del}}$, $X_f^{\text{cont}}$, and $X_f^{\text{del}}$ similarly. We have the inequality $0 \leq \rho_e + \rho_f - m \leq a-1$. It must be that $|X_e^{\text{del}}| \geq k-\rho_e$ and $|X_f^{\text{del}}| \geq k-\rho_f$; otherwise, $M$ would contain a proper $(U_{a_0,b_0})$-minor such that $a_0 < a$ and $b_0 > a_0$, contradiction. Let $Y_e$ be any subset of $X_e^{\text{del}}$ containing exactly $k-\rho_e$ elements and let $Y_f$ be any subset of $X_f^{\text{del}}$ containing exactly $k-\rho_f$ elements. Let $N = \nat \backslash (Y_e \cup Y_f)$. Note that $N$ contains $M$ as a minor. The rank of $N$ is $m$, and $|E(N)| = \rho_e + \rho_f$. This implies $\eta(N) = \rho_e + \rho_f - m$, which is strictly less than $a$ by assumption, contradiction.
        \item[(2)] This implies $0 \leq k+\rho_e-m \leq a-1$. It must be that $|X_e^{\text{del}}| \geq k-\rho_e$; otherwise, $M$ would contain a proper $(U_{a_0,b_0})$-minor such that $a_0 < a$ and $b_0 > a_0$, contradiction. Let $Y_e$ be any subset of $X_e^{\text{del}}$ containing exactly $k-\rho_e$ elements. Let $N = \nat \backslash Y_e$. Note that $N$ contains $M$ as a minor. The rank of $N$ is $m$, and $|E(N)| = k+\rho_e$. This implies $\eta(N) = k+\rho_e-m$, which is strictly less than $a$ by assumption, contradiction. 
        \item[(3), (5)] Let $X_e^{\text{cont}}$ and $X_e^{\text{del}}$ be disjoint subsets of $X_e$ with sizes $m-\rho_f$ and $k-\rho_e$ respectively. The matroid $M = S(\nat )/X_e^{\text{cont}}\backslash X_e^{\text{del}}$ consists of $k-m+\rho_e+\rho_f$ clones in rank $\rho_f$. We have
    \begin{align*}
        \rho_f &\geq (k+1)-a \geq b-a, \text{ and}\\
        k-m+\rho_e+\rho_f & \geq k+ (a-k) + (b-a) = b.
    \end{align*}
        The nullity $\eta(M)$ is equal to $k-m+\rho_e$ which is greater than or equal to $a$ as $m-\rho_e \in [a,k-a]$. Therefore, $M$ contains $\dunif $ as a proper restriction. \qedhere
    \end{enumerate}
\end{proof}

\begin{corollary}\label{cor-{num of doubleton excluded minors}}\leavevmode
Let $k \geq 2(b-a)$. The total number of doubleton excluded minors for $\class{k}$ is
\begin{align}\label{eq-{total number doubleton excluded minors}}
    \frac{1}{6}a(-2a^2 + 3ak + 3k + 2).
\end{align}
\end{corollary}
\begin{proof} The number of excluded minors from row (3) of \ref{tab-{doubleton excluded minors for unif}} is
\begin{align}\label{eq-{row 3 excluded minors}}
    \sum\limits_{\rho_e = 1}^{a-1}\left(\sum\limits_{i = 1}^{\rho_e}i\right) = \frac{1}{6}(a-1)(a)(a+1)
\end{align}
and the number of excluded minors from row (5) is
\begin{align}\label{eq-{row 5 excluded minors}}
    \sum\limits_{\rho_e = k-a+1}^{k}\left(\sum\limits_{i = 1}^{k-\rho_e + 1}(k-a+1)\right) = \frac{1}{2}a(-a^2 + ak+k+1).
\end{align}
The expression in \ref{eq-{total number doubleton excluded minors}} is equal to the sum of the quantities in Equations \ref{eq-{row 3 excluded minors}} and \ref{eq-{row 5 excluded minors}}.
\end{proof}

If $\rho = \rho_{(\rho_e,\rho_f)}^m \in Ex\big(\class{k}\big)$, we denote $\rho$ as $Ex_{(\rho_e,\rho_f)}^m$.

\begin{example}\label{ex-{U37}}
    Let $\cal{C}$ be the class of matroids characterized by $Ex(\cal{C}) = \{U_{3,7}, U_{4,7}\}$. Let $k =8$. The singleton excluded minors for $\class{8}$ are:
    \begin{itemize}
        \item $Ex^3$, $Ex^4$, $Ex^5$
    \end{itemize}
    The doubleton excluded minors for $\class{8}$ are:
    \begin{itemize}
        \item $Ex_{(1,6)}^6$, $Ex_{(2,6)}^6$, $Ex_{(2,6)}^7$, $Ex_{(2,7)}^7$
        \item $Ex_{(6,6)}^m$ for $m \in [6,11]$
        \item $Ex_{(6,7)}^m$ and $Ex_{(7,7)}^m$ for $m \in [7,12]$
        \item $Ex_{(6,8)}^m$, $Ex_{(7,8)}^m$, and $Ex_{(8,8)}^m$ for $m \in [8,13]$
    \end{itemize}
\end{example}

\begin{theorem}\label{thm-{bound on complexity of uniform}}
    For $b \geq 2a$, let $\cal{C}$ be the class of matroids characterized by
    \begin{align*}
        Ex(\cal{C}) = \{\unif , \dunif \}.
    \end{align*}
   Let $k \geq 2(b-a)$ and let $\rho \in Ex\big(\class{k}\big)$. Then $\rho$ has an $(a-1)$-corner decomposition.
\end{theorem}

\begin{proof} 
Let $(\{e\}, \rho)$ be a $k$-polymatroid in $\class{k}$. By Lemma \ref{prop-{general unif excluded singletons}}, $\rho(\{e\}) \leq a-1$ or $k-\rho(\{e\}) \leq a-1$. We claim that $\rho$ has an $(a-1)$-corner decomposition 
\begin{align*}
    \rho = \tau + (k-(a-1))r.
\end{align*}
The following table tracks $\tau$ and $r$ for each choice of $\rho(\{e\})$:
\begin{table}[H]
\centering
\begin{tabular}{ |c|c|c| }
 \hline
$\rho(\{e\})$ value (or range)& $\tau$&$r$\\
\hline
$0$&$U_{0,1}$&$U_{0,1}$\\
 \hline
 $[1, a-1]$&$\rho(\{e\})\coloop$&$U_{0,1}$\\
 \hline
 $k-(a-1)$&$U_{0,1}$&$\coloop$\\
 \hline
 $[k-a+2,k]$&$(\rho(\{e\})-k+a-1)\coloop$&$\coloop$\\
 \hline
\end{tabular}
\end{table}
It is easy to check that $\tau$ is an $(a-1)$-polymatroid and $r$ is a maximally-separated matroid. In particular, for the final row, since $\rho(\{e\}) \in [0,k]$, it must be that $\rho(\{e\})-k+a-1\leq a-1$. 

Next, consider $|E| = 2$ where $E = \{e, f\}$. Let $\betaInv := \rho_e+\rho_f - m$.\footnote{The \textit{signed beta invariant} of a matroid $M = (E,r)$ is $\betaInv(M) = -\sum\limits_{X \subseteq E}(-1)^{|X|}r(X)$. \cite{ardila benedetti doker}} For each choice of $\rho$ satisfying the hypothesis of Theorem \ref{thm-{bound on complexity of uniform}}, we can write $\rho = \tau + (k-(a-1))r$ where $r$ is given in the table below, and
\begin{align*}
    \tau = \betaInv U_{1,2} + \left(\left(\rho_e - \betaInv\right)\coloop \oplus \left(\rho_f - \betaInv \right)\coloop\right).
\end{align*}
(Note: Direct sums are written in the form $\rho|_{\{e\}}\oplus \rho|_{\{f\}}$.) We use Proposition \ref{prop-{general unif excluded doubletons}} to determine the allowed possibilities for the first three columns. Any proper minor is for each allowable $\rho$ is either empty or falls into the $|E| = 1$ case so it suffices to address $\rho$ itself.

\begin{table}[H]
\centering
\begin{tabular}{ |c|c|c|c| }
\hline
$\rho_e$ range & $\rho_f$ range & $m$ range & $r$\\
\hline
$[0,a-1]$ & $[0,a-1]$&$[\rho_f, \rho_e+\rho_f]$&$U_{0,1}  \oplus U_{0,1} $\\
\hline
$[0,a-1]$& $[k-a+1, k]$&$[\rho_e + (k-a+1), \rho_e+\rho_f]$&$U_{0,1}  \oplus \coloop$\\
\hline
$[k-a+1,k]$ & $[k-a+1,k]$&$[\rho_f + (k-a+1), \rho_e + \rho_f]$&$\coloop \oplus \coloop$\\
\hline
\end{tabular}
\end{table}

    We look to Remark \ref{rem-{sequence of compressions}} and Lemma \ref{lem-{Gamma implies small E}}. Let $\Gamma_b$ be the set of polymatroids in $Ex(\class{k})$ on a ground set of size $b$. Assume there is some $\rho \in \class{k}$ which is a decompression of a polymatroid in $\Gamma_b$; then $|E(\rho)| = b+1$. Since $\rho$ is an excluded minor, it cannot properly contain a singleton or doubleton minor isomorphic to $(\unif )$ or $(\dunif )$. Therefore, if $\rho'$ is a singleton or doubleton minor of $\rho$, then $\rho'$ must be isomorphic to one of the options given in the tables above. Each of those options has an $(a-1)$-corner decomposition as shown. By Proposition \ref{prop-{singleton and doubleton minors}}, $\rho$ has an $(a-1)$-corner decomposition.
\end{proof}

\begin{example}
Let $\cal{C}$ be the class of matroids from Example \ref{ex-{U37}}. Let $\rho$ be a singleton $k$-polymatroid in $\class{8}$ with rank $\rho_e$. If $\rho_e \in \{0,1,2\}$, then $r = U_{0,1} $ (loop) and $\tau = \rho_e\coloop$, and $\rho$ is confined to the $2$-corner given by $[[0,2]]$. If $\rho_e \in \{6, 7, 8\}$, then $r = \coloop$ (coloop) and $\tau = (\rho_e-6)\coloop$, and $\rho$ is confined to the $2$-corner given by $[[6,8]]$.

$$\begin{tikzpicture}
		[ lightdashed/.style={very thin,gray, dashed},
            darkdashed/.style={very thick, black, dashed},
            darksolid/.style = {very thick, black},
            lightsolid/.style = {very thin, gray},
			axis/.style={->,black, very thin}]
			
\draw[axis] (0,0)--(8.5,0) node[anchor=west]{$e$};

\node[shape = circle, fill = black, scale = 0.5] at (0,0)  {}; 
    \node[] at (0,0.3)  {$\rho_e = 0$}; 
\node[shape = circle, fill = black, scale = 0.5] at (1,0)  {}; 
    \node[] at (1,0.3)  {$1$}; 
\node[shape = circle, fill = black, scale = 0.5] at (2,0)  {}; 
    \node[] at (2,0.3)  {$2$}; 
\node[shape = circle, fill = black, scale = 0.5] at (6,0)  {}; 
    \node[] at (6,0.3)  {$6$}; 
\node[shape = circle, fill = black, scale = 0.5] at (7,0)  {}; 
    \node[] at (7,0.3)  {$7$}; 
\node[shape = circle, fill = black, scale = 0.5] at (8,0)  {}; 
    \node[] at (8,0.3)  {$8$}; 
    
\end{tikzpicture}$$

Now let $E = \{e,f\}$. Below is a complete list\footnote{With some duplicates -- an entry is enclosed in square brackets if and only if an isomorphic polymatroid appears earlier (i.e. to the left of or above it) in the table.} of the $8$-polymatroids $\rho = \rho_{(\rho_e,\rho_f)}^m$ in $\class{8}$ up to isomorphism. The entry in the $i$th row and $j$th column corresponds to $\rho = \tau + 7r$ where the $2$-polymatroid $\tau$ is given in the $i$th row of the first column and the maximally-separated matroid $r$ is given in the $j$th column of the first row. Under each $r$, the $2$-corner to which $B_\rho$ is confined is given in parentheses. Note: Any $\rho \in \class{8}$ confined to the $2$-corner given by $[[6,8]] \times [[0,2]]$ is isomorphic to a polymatroid in column two. 

\begin{table}[H]
\centering
\begin{tabular}{|c|c|c|c|c|}
\hline
\backslashbox{Row}{Col}&$0$& $1$&  $2$&  $3$\\
\hline
$0$&\backslashbox{$\tau$}{$\substack{r\\ \text{($2$-corner)}}$}
&$\substack{U_{0,1}  \oplus U_{0,1} \\ ([0,2] \times [0,2])}$ & $\substack{U_{0,1}  \oplus \coloop\\ \left([0,2] \times [6,8]\right)}$ & $\substack{\coloop \oplus \coloop\\ \left([6,8] \times [6,8]\right)}$\\
\hline
 $1$ & $U_{0,1}  \oplus U_{0,1} $&$\rho_{(0,0)}^0$&$\rho_{(0,6)}^6$&$\rho_{(6,6)}^{12}$\\
\hline
 $2$ &$U_{0,1}  \oplus \coloop$ &$\rho_{(0,1)}^1$&$\rho_{(0,7)}^7$&$\rho_{(6,7)}^{13}$\\
\hline
 $3$ &$U_{0,1}  \oplus 2\coloop$ &$\rho_{(0,2)}^2$&$\rho_{(0,8)}^8$&$\rho_{(6,8)}^{14}$\\
\hline
 $4$ &$\coloop \oplus U_{0,1} $ &$[\rho_{(1,0)}^1]$&$\rho_{(1,6)}^7$&$[\rho_{(7,6)}^{13}]$\\
\hline
 $5$ &$\coloop \oplus \coloop$ &$\rho_{(1,1)}^2$&$\rho_{(1,7)}^8$&$\rho_{(7,7)}^{14}$\\
\hline
 $6$ &$\coloop \oplus 2\coloop$ &$\rho_{(1,2)}^3$&$\rho_{(1,8)}^9$&$\rho_{(7,8)}^{15}$\\
\hline
 $7$ &$2\coloop \oplus U_{0,1} $&$[\rho_{(2,0)}^2]$ &$\rho_{(2,6)}^8$&$[\rho_{(8,6)}^{14}]$ \\
\hline
 $8$ &$2\coloop \oplus \coloop$&$[\rho_{(2,1)}^3]$ &$\rho_{(2,7)}^9$&$[\rho_{(8,7)}^{15}]$ \\
\hline
 $9$ &$2\coloop \oplus 2\coloop$&$\rho_{(2,2)}^4$&$\rho_{(2,8)}^{10}$&$\rho_{(8,8)}^{16}$\\
\hline
 $10$ &$U_{1,2}$ &$\rho_{(1,1)}^1$&$\rho_{(1,7)}^7$&$\rho_{(7,7)}^{13}$\\
\hline
 $11$ &$U_{1,2} + (U_{0,1}  \oplus \coloop)$ &$\rho_{(1,2)}^2$&$\rho_{(1,8)}^8$&$\rho_{(7,8)}^{14}$\\
\hline
 $12$ &$U_{1,2} + (\coloop \oplus U_{0,1} )$ &$[\rho_{(2,1)}^2]$ &$\rho_{(2,7)}^8$&$[\rho_{(8,7)}^{14}]$ \\
\hline
 $13$ &$U_{1,2} + (\coloop \oplus \coloop)$ &$\rho_{(2,2)}^3$&$\rho_{(2,8)}^9$&$\rho_{(8,8)}^{15}$\\
\hline
 $14$ &$2U_{1,2}$ &$\rho_{(2,2)}^2$&$\rho_{(2,8)}^8$&$\rho_{(8,8)}^{14}$\\
\hline
\end{tabular} 
\end{table} 

Every polymatroid $\rho$ in the table has base polytope $B_\rho$ confined to one of the white regions (the $2$-corners of $[[0,8]]\times [[0,8]]$) below.
$$
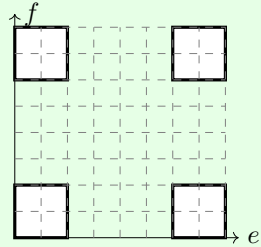
\begin{tikzpicture}
		[scale = 0.35, lightdashed/.style={very thin,gray, dashed},
            darkdashed/.style={very thick, black, dashed},
            darksolid/.style = {very thick, black},
            lightsolid/.style = {very thin, gray},
			axis/.style={->,black, very thin}]
			
\draw[axis] (0,0)--(8.5,0) node[anchor=west]{$e$};
\draw[axis] (0,0) -- (0,8.5) node[anchor=west]{$f$};

\draw[darksolid, fill = white] (0,0)--(2,0)--(2,2)--(0,2)--cycle;
\draw[darksolid, fill = white] (0,6)--(2,6)--(2,8)--(0,8)--cycle;
\draw[darksolid, fill = white] (6,0)--(8,0)--(8,2)--(6,2)--cycle;
\draw[darksolid, fill = white] (6,6)--(8,6)--(8,8)--(6,8)--cycle;

    \draw[lightdashed] (1,0)--(1,8);
    \draw[lightdashed] (2,0)--(2,8);
    \draw[lightdashed] (3,0)--(3,8);
    \draw[lightdashed] (4,0)--(4,8);
    \draw[lightdashed] (5,0)--(5,8);
    \draw[lightdashed] (6,0)--(6,8);
    \draw[lightdashed] (7,0)--(7,8);
    \draw[lightdashed] (8,0)--(8,8);
    \draw[lightdashed] (0,1)--(8,1);
    \draw[lightdashed] (0,2)--(8,2);
    \draw[lightdashed] (0,3)--(8,3);
    \draw[lightdashed] (0,4)--(8,4);
    \draw[lightdashed] (0,5)--(8,5);
    \draw[lightdashed] (0,6)--(8,6);
    \draw[lightdashed] (0,7)--(8,7);
    \draw[lightdashed] (0,8)--(8,8);

 \end{tikzpicture}$$

Depicted below are the polymatroids $\rho$ whose base polytopes $B_\rho$ lie in the $2$-corner given by $[[0,2]]\times [[0,2]]$. Each $B_\rho$ is labeled with $\tau_i$ where $i$ corresponds to the row of the table in which $\tau$ is found.
$$
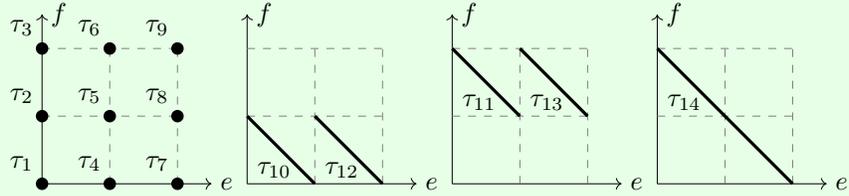
\begin{tikzpicture}
		[scale = 0.9, lightdashed/.style={very thin,gray, dashed},
            darkdashed/.style={very thick, black, dashed},
            darksolid/.style = {very thick, black},
            lightsolid/.style = {very thin, gray},
			axis/.style={->,black, very thin}]
			
\draw[axis] (0,0)--(2.5,0) node[anchor=west]{$e$};
\draw[axis] (0,0) -- (0,2.5) node[anchor=west]{$f$};

    \draw[lightdashed] (1,0)--(1,2);
    \draw[lightdashed] (2,0)--(2,2);
    \draw[lightdashed] (0,1)--(2,1);
    \draw[lightdashed] (0,2)--(2,2);

    \node[shape = circle, fill = black, scale = 0.5] at (0,0)  {}; 
    \node[] at (-0.3,0.3)  {$\tau_1$}; 
    \node[shape = circle, fill = black, scale = 0.5] at (0,1)  {}; 
    \node[] at (-0.3,1.3)  {$\tau_2$}; 
    \node[shape = circle, fill = black, scale = 0.5] at (0,2)  {};
    \node[] at (-0.3,2.3)  {$\tau_3$}; 
    \node[shape = circle, fill = black, scale = 0.5] at (1,0)  {}; 
    \node[] at (0.7,0.3)  {$\tau_4$}; 
    \node[shape = circle, fill = black, scale = 0.5] at (1,1)  {}; 
    \node[] at (0.7,1.3)  {$\tau_5$}; 
    \node[shape = circle, fill = black, scale = 0.5] at (1,2)  {}; 
    \node[] at (0.7,2.3)  {$\tau_6$}; 
    \node[shape = circle, fill = black, scale = 0.5] at (2,0)  {}; 
    \node[] at (1.7,0.3)  {$\tau_7$}; 
    \node[shape = circle, fill = black, scale = 0.5] at (2,1)  {}; 
    \node[] at (1.7,1.3)  {$\tau_8$}; 
    \node[shape = circle, fill = black, scale = 0.5] at (2,2)  {}; 
    \node[] at (1.7,2.3)  {$\tau_9$}; 

 \end{tikzpicture}\begin{tikzpicture}
		[scale = 0.9, lightdashed/.style={very thin,gray, dashed},
            darkdashed/.style={very thick, black, dashed},
            darksolid/.style = {very thick, black},
            lightsolid/.style = {very thin, gray},
			axis/.style={->,black, very thin}]
			
\draw[axis] (0,0)--(2.5,0) node[anchor=west]{$e$};
\draw[axis] (0,0) -- (0,2.5) node[anchor=west]{$f$};

    \draw[lightdashed] (1,0)--(1,2);
    \draw[lightdashed] (2,0)--(2,2);
    \draw[lightdashed] (0,1)--(2,1);
    \draw[lightdashed] (0,2)--(2,2);

    \draw[darksolid] (1,0)--(0,1);
    \node[] at (0.4,0.2)  {$\tau_{10}$}; 
    \draw[darksolid] (2,0)--(1,1);
    \node[] at (1.4,0.2)  {$\tau_{12}$}; 

 \end{tikzpicture}\begin{tikzpicture}
		[scale = 0.9, lightdashed/.style={very thin,gray, dashed},
            darkdashed/.style={very thick, black, dashed},
            darksolid/.style = {very thick, black},
            lightsolid/.style = {very thin, gray},
			axis/.style={->,black, very thin}]
			
\draw[axis] (0,0)--(2.5,0) node[anchor=west]{$e$};
\draw[axis] (0,0) -- (0,2.5) node[anchor=west]{$f$};

    \draw[lightdashed] (1,0)--(1,2);
    \draw[lightdashed] (2,0)--(2,2);
    \draw[lightdashed] (0,1)--(2,1);
    \draw[lightdashed] (0,2)--(2,2);

    \draw[darksolid] (1,1)--(0,2);
    \node[] at (0.4,1.2)  {$\tau_{11}$}; 
    \draw[darksolid] (2,1)--(1,2);
    \node[] at (1.4,1.2)  {$\tau_{13}$}; 

 \end{tikzpicture}\begin{tikzpicture}
		[scale = 0.9, lightdashed/.style={very thin,gray, dashed},
            darkdashed/.style={very thick, black, dashed},
            darksolid/.style = {very thick, black},
            lightsolid/.style = {very thin, gray},
			axis/.style={->,black, very thin}]
			
\draw[axis] (0,0)--(2.5,0) node[anchor=west]{$e$};
\draw[axis] (0,0) -- (0,2.5) node[anchor=west]{$f$};

    \draw[lightdashed] (1,0)--(1,2);
    \draw[lightdashed] (2,0)--(2,2);
    \draw[lightdashed] (0,1)--(2,1);
    \draw[lightdashed] (0,2)--(2,2);

    \draw[darksolid] (2,0)--(0,2);
    \node[] at (0.4,1.2)  {$\tau_{14}$}; 

 \end{tikzpicture}$$

Every base polytope from one of the other three $2$-corners of $[[0,8]] \times [[0,8]]$ is equivalent to one of the $\tau_i$s above. 

\end{example}

\end{document}